\documentclass{article}

\usepackage[preprint,nonatbib]{neurips_2020}

\usepackage{appendix}
\usepackage{algorithmic}
\usepackage{url}            
\usepackage{amsfonts}       
\usepackage{nicefrac}       
\usepackage{microtype}      
\usepackage{microtype}
\usepackage{mathrsfs}
\usepackage{graphicx}
\usepackage{subfigure}
\usepackage{amssymb}
\usepackage{booktabs} 
\usepackage{bm}
\usepackage{url}            
\usepackage{amsfonts}       
\usepackage{nicefrac}       
\usepackage{microtype}      
\usepackage{bm}
\usepackage{amsmath}
\usepackage{amssymb}
\usepackage{mathtools}
\usepackage{titlesec}
\usepackage{xcolor}
\usepackage{wrapfig}
\usepackage{caption}
\usepackage{cite}
\usepackage{enumitem}
\usepackage[utf8]{inputenc} 
\usepackage[T1]{fontenc}
\usepackage{setspace}
\usepackage{multirow}
\usepackage{comment}
\usepackage[colorlinks=true,bookmarksopen,pdfsubject={algorithms},linkcolor={blue}, anchorcolor={black}, citecolor={blue}, filecolor={magenta}, menucolor={black},backref=none,urlcolor={blue}]{hyperref}

\newtheorem{thm}{Theorem}

\newtheorem{proposition}{Proposition}
\newtheorem{corollary}{Corollary}

\title{A Novel Evolution Strategy with Directional Gaussian Smoothing for Black-box Optimization}

\author{%
Jiaxin Zhang\\
  \small Computer Science and Mathematics Division\\
  \small Oak Ridge National Laboratory\\
  \small \texttt{zhangj@ornl.gov} 
  \And
Hoang Tran\\
  \small Computer Science and Mathematics Division\\
  \small Oak Ridge National Laboratory\\
  \small \texttt{tranha@ornl.gov}   
  \And
 Dan Lu\\
 \small  Computational Science and Engineering Division\\
 \small  Oak Ridge National Laboratory\\
  \small \texttt{lud1@ornl.gov} \\
  \And
Guannan Zhang\thanks{Corresponding Author}\\
  \small Computer Science and Mathematics Division\\
  \small Oak Ridge National Laboratory\\
  \small \texttt{zhangg@ornl.gov}   
 }
\begin{document}

\maketitle

\begin{abstract}
We propose an improved evolution strategy (ES) using a novel nonlocal gradient operator for high-dimensional black-box optimization. Standard ES methods with $d$-dimensional Gaussian smoothing suffer from the curse of dimensionality due to the high variance of Monte Carlo (MC) based gradient estimators. To control the variance, Gaussian smoothing is usually limited in a small region, so existing ES methods lack nonlocal exploration ability required for escaping from local minima. We develop a nonlocal gradient operator with directional Gaussian smoothing (DGS) to address this challenge. The DGS conducts 1D nonlocal explorations along $d$ orthogonal directions in $\mathbb{R}^d$, each of which defines a nonlocal directional derivative as a 1D integral. We then use Gauss-Hermite quadrature, instead of MC sampling, to estimate the $d$ 1D integrals to ensure high accuracy (i.e., small variance). Our method enables effective nonlocal exploration to facilitate the global search in high-dimensional optimization. We demonstrate the superior performance of our method in three sets of examples, including benchmark functions for global optimization, and real-world science and engineering applications.
\end{abstract}

\section{Introduction}
\label{sec:intro}

Evolution strategy (ES) is a type of evolutionary algorithms for black-box optimization, where we search for the optima of a $d$-dimensional loss function $F(\bm x)$ given access to only its function queries. This is motivated by several applications where the loss function's gradient is inaccessible, e.g., in optimizing neural network architecture \cite{Real_ICML17,Miikkulainen_ICML17}, reinforcement learning \cite{Houthooft18,Ha_Schmidhuber,SHCS17}, and design of adversarial attacks to deep networks \cite{10.1145/3128572.3140448}. There are several types of evolutionary algorithms, including 
genetic algorithms \cite{10.5555/534133}, differential evolution \cite{DAS20161}, natural evolution strategies \cite{Wierstra_NES_14}, neuroevolution \cite{Stanley_neuroevolution}, and covariance matrix adaptation ES (CMA-ES) 
\cite{Hansen_Ostermeier_CMA_01,Hansen_CMA, 1554902,8410043, 10.1145/2908812.2908863}.

In this work, we consider a particular class of ES that is based on Gaussian smoothing (GS) (e.g., \cite{SHCS17, Mania2018SimpleRS, Maheswaranathan_GuidedES, Choromanski_ES-Active}). 
GS-based ES first smooths the landscape of the loss function with $d$-dimensional Gaussian convolution and then estimates the gradient of the {smoothed} loss function using a random population generated by Monte Carlo (MC) sampling. 
GS-based ES is a promising strategy to handle loss landscapes possessing many local minima. Theoretically, the Gaussian convolution with a large smoothing radius enables {nonlocal} exploration that reduces the local minima effect and improves global structure characterization. However, practically, such nonlocal exploration is limited in low- to medium-dimensional settings because of the high variance (i.e., low accuracy) of MC estimation. When GS has a large smoothing radius in a high-dimensional space, the exploration domain (i.e., the high probability region) will be significantly large. To reduce the variance thus ensure the estimation accuracy, MC sampling requires either a prohibitively large number of function evaluations or a very small smoothing radius. Given the computing constraint of the former, the latter is usually applied, but this results in local exploration. Therefore, GS-based ES has been mostly used to estimate {\em local} gradients in high-dimensional black-box optimization.  

Several studies have been performed to improve the gradient estimation in GS-based ES. Most of them center on enhancing MC estimators, 
such as by variance reduction strategies \cite{MWDS18,CRSTW18,10.5555/3326943.3326962}, exploiting historical data \cite{Maheswaranathan_GuidedES,Meier_OPTRL_2019}, employing active subspaces \cite{Choromanski_ES-Active}, and searching on latent low-dimensional manifolds \cite{Sener2020Learning}. Despite of some improvements, these techniques did not fundamentally solve the limitations in GS-based ES caused by the MC estimation. 

We develop a novel nonlocal gradient operator based on 
{directional Gaussian smoothing} (DGS), and use it to replace the standard GS operator in ES to enhance nonlocal exploration in high dimensions. We name our new operator as DGS gradient and the DGS-based ES algorithm as DGS-ES method.
The key idea of the DGS gradient is to conduct 1D nonlocal explorations along $d$ orthogonal directions in $\mathbb{R}^d$, each of which defines a nonlocal directional derivative as a 1D integral. Then we use Gauss-Hermite quadrature, instead of MC sampling, to estimate the $d$ 1D integrals to ensure high accuracy. Next, the estimated directional derivatives are assembled to form the DGS gradient. Compared with existing methods, our DGS-ES approach can achieve long-range exploration by being able to use a large smoothing radius and meanwhile obtains high estimation accuracy of gradients through accurate integral calculation. 

The proposed DGS gradient is a new {operator} for identifying search directions, not an estimator of the local gradient. In the local setting (the smoothing radius approaching to zero), we verified that the DGS-gradient operator is consistent with the GS-based gradient. However, in the nonlocal setting, the DGS gradient {is significantly different from the GS-based gradient in that it is feasible to obtain an accurate estimator of the DGS gradient with large smoothing radius} for nonlocal exploration.

\textbf{Summary of contributions.} Our contribution in this paper is three fold: (1) we develop the DGS gradient operator for effective nonlocal exploration in ES, which advances the global search in high-dimensional black-box optimization; (2) we develop an accurate estimator for the DGS gradient using Gauss-Hermite quadrature, 
%
which accelerates the convergence of ES; and (3) we demonstrate the superior performance of our method on both high-dimensional, non-convex benchmark optimization problems, and two real-world science and engineering applications.

\textbf{Related works.}
The literature on black-box optimization is extensive. We review three types of methods that are closely related to this work (see \cite{RiosSahinidis13,Larson_et_al_19} for thorough reviews).
%
%
(1) \textit{Random search}. This type of methods randomly generate the next search direction and estimate the directional derivative {or perform direct search for the updates}. Examples are two-point approaches \cite{FKM05,NesterovSpokoiny15,Duchi2015OptimalRF, MAL-024}, 
coordinate-descent algorithms \cite{10.5555/2999325.2999433}, 
three-point methods \cite{Bergou2019}, and binary search with adaptive radius \cite{Golovin2020GradientlessDH}. From theoretical perspective, { an analysis of two-point schemes based on GS is presented in the seminal paper \cite{NesterovSpokoiny15}, extended in \cite{doi:10.1137/120880811} for non-convex and in \cite{10.5555/3122009.3153008} for non-smooth loss functions.} Existing studies also focus on estimating local derivatives rather than nonlocal exploration. 
%
%
(2) {\em Local gradient estimation}. The most straightforward way is to use finite differences. An alternative is to use linear interpolation in a small neighborhood of the current state to estimate local gradients \cite{BCCS19b}. Another way is to estimate the local gradient by averaging multiple directional estimates by two-point schemes, and the GS-based ES methods \cite{SHCS17,Mania2018SimpleRS, Maheswaranathan_GuidedES,Choromanski_ES-Active} can be assigned to this type of methods.
It is possible to augment ES by integrating the estimated gradient with new gradient-based algorithms, such as ADMM \cite{2017arXiv171007804L}, adaptive momentum method \cite{Chen2019ZOAdaMMZA}, and conditional gradient \cite{10.5555/3327144.3327264}. A comparison of local gradient estimation methods can be found in \cite{BCCS19}. 
%
%
(3) \textit{Smoothing techniques}. Sphere smoothing is a method similar to GS and was discussed in \cite{FKM05}. Analysis of GS applied to step functions are presented in \cite{doi:10.1080/10556780500140029}. Other strategies to transform a nonconvex and noisy optimization to a convex or more friendly version are $p$-th power transformation \cite{Li-convexification} and $\ell^2$ regularization \cite{Carlsson19}. An algorithm for estimating computational noise affecting a smooth simulation is developed in \cite{Mor2011EstimatingCN}.



%

\section{Black-box optimization}\label{sec:setting}
We are interested in solving the following black-box optimization problem
\begin{equation}\label{e1}
    \min_{\bm x \in \mathbb{R}^d} F(\bm x),
\end{equation}
where $\bm x = (x_1, \ldots, x_d) \in \mathbb{R}^d$ consists of $d$ tuning parameters, and $F: \mathbb{R}^d \rightarrow \mathbb{R}$ is a $d$-dimensional black-box loss function.
We assume that the gradient $\nabla F(\bm x)$ is unavailable, and $F(\bm x)$ is only accessible via function evaluations. 

We briefly recall the class of ES methods \cite{SHCS17} that use GS \cite{FKM05,NesterovSpokoiny15} to estimate local gradients. The smoothed loss is defined by
$
       F_{\sigma}(\bm x) = \mathbb{E}_{\bm u \sim \mathcal{N}(0, \mathbf{I}_d)} \left[F(\bm x + \sigma \bm u) \right], 
$
where $\mathcal{N}(0, \mathbf{I}_d)$ is the $d$-dimensional standard Gaussian distribution, and $\sigma > 0$ is the smoothing radius. $F_{\sigma}(\bm x)$ inherits many characteristics from $F(\bm x)$, e.g., convexity, the Lipschitz constant. Moreover, for any $\sigma >0$, $F_{\sigma}$ is always differentiable even if $F$ is not. 
%
%
%
%
%
The standard ES method \cite{SHCS17} represents the $\nabla F_\sigma(\bm x)$ as an expectation and estimate it by drawing $M$ random samples $\{\bm u_m\}_{m=1}^M$
from $\mathcal{N}(0,\mathbf{I}_d)$, i.e., 
\begin{equation}\label{e40}
    \nabla F_{\sigma}(\bm x) = 
    \frac{1}{\sigma}\mathbb{E}_{\bm u \sim \mathcal{N}(0, \mathbf{I}_d)} \left[F(\bm x + \sigma \bm u)\, \bm u\right] \approx \frac{1}{M\sigma}\sum_{m=1}^M F(\bm x + \sigma \bm u_m)\bm u_m.
\end{equation}
%
Then, the MC estimator is substituted into any gradient-based algorithm to update the state $\bm x$.

%

\section{New method: an evolution strategy with directional Gaussian smoothing}
\label{sec:DGS-ES}
We present our main contributions in this section. We start by introducing the DGS gradient in \S \ref{sec:grad}. 
%
In \S \ref{sec:ada_DGS-ES}, we introduce an accurate estimator of the DGS gradient and describe the proposed DGS-ES algorithm in detail. 
The proposed method was inspired by the following key idea:

{\em {\bf Key idea}: The DGS gradient conducts 1D nonlocal explorations along $d$ orthogonal directions in $\mathbb{R}^d$, each of which defines a nonlocal directional derivative as a 1D integral. The Gauss-Hermite quadrature, instead of MC sampling, is used to estimate the $d$ 1D integrals to achieve high accuracy.}

\subsection{The nonlocal DGS gradient operator}\label{sec:grad}
To proceed, we first define a \emph{one-dimensional} cross section of $F(\bm x)$ as
\begin{equation*}
G(y \,| \,{\bm x, \bm \xi}) = F(\bm x + y\, \bm \xi), \;\; y \in \mathbb{R},
\end{equation*}
where $\bm x$ is the current state of $F(\bm x)$ and $\bm \xi$ is a unit vector in $\mathbb{R}^d$. Note that $\bm x$ and $\bm \xi$ can be viewed as parameters of the function $G$. We define the Gaussian smoothing of $G(y)$, denoted by $G_\sigma(y)$, by
\begin{equation}
\label{eq10}
    G_{\sigma}(y \,| \,{\bm x, \bm \xi}) :=  \frac{1}{\sqrt{2\pi}} \int_{\mathbb{R}} G(y + \sigma v\, |\, \bm x, \bm \xi)\, {\rm e}^{-\frac{v^2}{2}}\, dv
     = \mathbb{E}_{v \sim \mathcal{N}(0, 1)} \left[G(y + \sigma v\, |\, \bm x, \bm \xi) \right],
\end{equation}
which is the Gaussian smoothing of $F(\bm x)$ along the direction $\bm \xi$ in the neighbourhood of $\bm x$. 
The derivative of $G_{\sigma}(y|\bm x,\bm \xi)$ at $y = 0$ can be represented by a one-dimensional expectation
\begin{equation}\label{e4}
    \mathscr{D}[G_{\sigma}(0 \,|\, \bm x, \bm \xi)] 
     = \frac{1}{\sigma}\,\mathbb{E}_{v \sim \mathcal{N}(0,1)} \left[G(\sigma v \, | \, \bm x, \bm \xi)\, v\right],
\end{equation}
where $\mathscr{D}[\cdot]$ denotes the differential operator. The difference between the directional derivative of $F_\sigma(\bm x)$ and Eq.~\eqref{e4} is that 
$\mathscr{D}[G_{\sigma}(0 \,|\, \bm x, \bm \xi)]$ 
only involves the directionally smoothed function in Eq.~\eqref{eq10}. 

For a matrix $\bm \Xi := (\bm \xi_1, \ldots, \bm \xi_d)$ consisting of 
$d$ orthonormal vectors, we can define $d$ directional derivatives like those in Eq.~\eqref{e4} and assemble our DGS gradient as 
%
%
\begin{equation}\label{dev_smooth_func}
 \hspace{-0.8cm} \text{\bf The DGS gradient:}\quad {\nabla}_{\sigma, \bm \Xi}[F](\bm x) := \Big[\mathscr{D}[G_{\sigma}(0 \, |\, \bm x, \bm \xi_1)], \cdots, {\mathscr{D}}[G_{\sigma}(0\, |\, \bm x, \bm \xi_d)]\Big]\, \bm \Xi,
\end{equation}
where the orthogonal system $\bm \Xi$ and the smoothing radius $\sigma$ can be adjusted during an optimization process. Next we describe how to integrate the DGS gradient into ES framework. 



\subsection{The DGS-ES algorithm}\label{sec:ada_DGS-ES}
The key step to integrate the DGS gradient in Eq.~\eqref{dev_smooth_func} into ES is to develop an accurate estimator. 
We exploit that each component of ${\nabla}_{\sigma, \bm \Xi}[F](\bm x)$ only involves a 1D integral, such that the Gauss-Hermite quadrature rule \cite{2013JSV...332.4403B,Handbook} can be used to approximate the integrals with high accuracy (shown in Eq.~\eqref{GH_error}). By doing a simple change of variable in Eq.~\eqref{e4}, the GH rule can be directly used to obtain the following estimator for $\mathscr{D}[G_{\sigma}(0 \,|\, \bm x, \bm \xi)]$, i.e.,
%
\begin{align}
 \widetilde{\mathscr{D}}^M[G_\sigma(0 \, | \, \bm x, \bm \xi)]  
   =
      \frac{1}{\sqrt{\pi}\sigma} \sum_{m = 1}^M w_m \,F(\bm x + \sqrt{2}\sigma v_m \bm \xi)\sqrt{2}v_m, \label{e8}
\end{align}
where $\{v_m\}_{m=1}^M$ are the roots of the $M$-th order Hermite polynomial
 and $\{w_m\}_{m=1}^M$ are quadrature weights.  
Both $v_m$ and $w_m$ can be found online\footnote{Nodes and weights for GH quadrature: \url{https://keisan.casio.com/exec/system/1281195844}} or in \cite{Handbook}. Compared with MC sampling, the error of Eq.~\eqref{e8} can be bounded by 
\begin{align}
\label{GH_error}
\hspace{-0.1cm}\big|(\widetilde{\mathscr{D}}^M- \mathscr{D})[G_\sigma] \big| \le C\frac{M\,!\sqrt{\pi}}{2^M(2M)\,!} \sigma^{2M-1}, 
\end{align}
where $M!$ is the factorial of $M$ and the constant $C>0$ is independent of $M$ and $\sigma$. 
Applying the GH quadrature rule $\widetilde{\mathscr{D}}^M$ to each component of ${\nabla}_{\sigma, \bm \Xi}[F](\bm x)$ in Eq.~\eqref{dev_smooth_func}, we define the following estimator: 
\begin{equation}\label{e5}
  \hspace{-0.2cm} \text{\bf The DGS estimator:}\quad \widetilde{\nabla}^M_{\sigma, \bm \Xi}[F](\bm x) = \Big[\widetilde{\mathscr{D}}^M[G_{\sigma}(0 \, |\, \bm x, \bm \xi_1)], \cdots, \widetilde{\mathscr{D}}^M[G_{\sigma}(0\, |\, \bm x, \bm \xi_d)]\Big]\, \bm \Xi.
\end{equation}
\begin{wrapfigure}{r}{0.5\textwidth}
    \vspace{-0.4cm}
   \hspace{-0.25cm} \includegraphics[scale = 0.45]{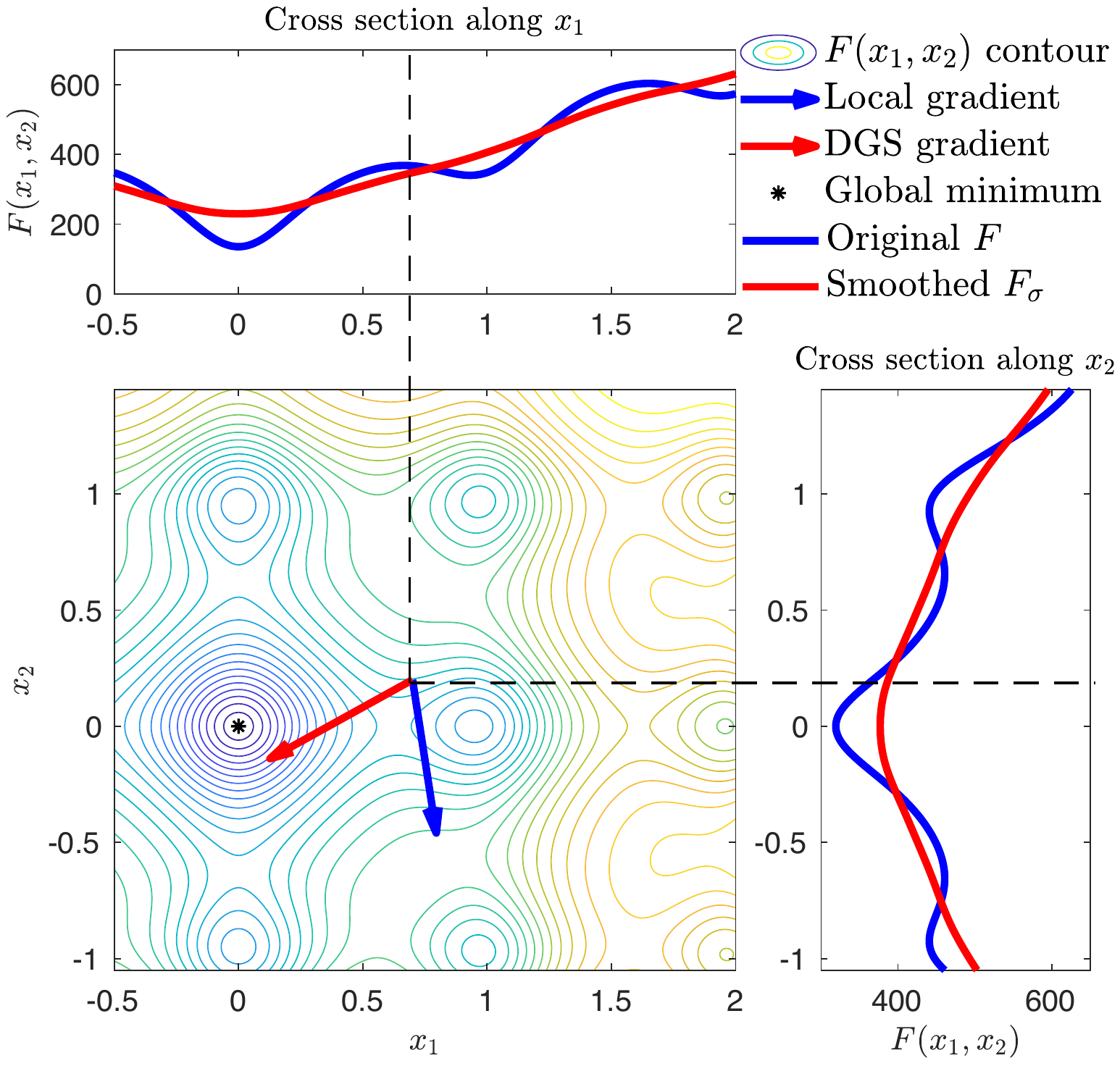}\vspace{-0.2cm}
    \caption{Illustration of the nonlocal exploration capability of our DGS gradient. In the central plot, the blue arrow points to the {local} gradient direction and the red arrow points to the DGS gradient direction. The top and right plots show the directionally smoothed functions along the two axes. Because the DGS gradient captures the nonlocal features of $F$, it can point to a direction much closer to the global minimum than the local gradient.}
    \label{fig0}
    \vspace{-1.4cm}
\end{wrapfigure}
The DGS estimator has the following features:
\vspace{-0.1cm}
\begin{itemize}[leftmargin=13pt]
    \item {\bf \em Nonlocality}: The directional smoothing allows for a large radius $\sigma$ to capture global structures of loss landscapes and help escape from local minima (illustrated in Figure \ref{fig0}).
    \item {\bf \em Accuracy}: The GH quadrature with the error bounded in Eq.~\eqref{GH_error} provides an estimator having much higher accuracy than MC, even when a large smoothing radius $\sigma$ is used.
    \item {\bf \em Portability}: The DGS gradient can be integrated into majority of gradient-based algorithms, e.g., gradient descent, Adam, and those with constraints (shown in \S\ref{sec:ex_3}).
    \item {\bf \em Scalability}: The DGS estimator in Eq.~\eqref{e5} requires $M\times d$ evaluations of $F(\bm x)$, and these evaluations are completely parallelizable as those in random sampling.
\end{itemize}



\paragraph{Random perturbation of $\bm \Xi$ and $\sigma$.}
The estimator in Eq.~\eqref{e5} is deterministic for a fixed $\bm \Xi$ and $\sigma$, making our approach short of random exploration.
%
%
%
To alleviate this issue, we add random perturbations to $\bm \Xi$ and $\sigma$. First, we add a small random rotation $\Delta \bm\Xi$ to $\bm \Xi$. To make $\bm \Xi + \Delta \bm\Xi$ orthonormal, i.e., 
$
   \mathbf{I}_d  = (\bm \Xi + \Delta \bm\Xi)^{\top} (\bm \Xi + \Delta \bm\Xi),
$ 
we generate $\Delta \bm \Xi$ as a random skew-symmetric matrix $\Delta \bm \Xi = -\Delta \bm \Xi^{\top}$ with small-value entries (controlled by $\alpha>0$), which will cancel out the first-order terms in $(\bm \Xi + \Delta \bm\Xi)^{\top} (\bm \Xi + \Delta \bm\Xi)$. The Gram-Schmidt operation is then used to eliminate the second-order term $\Delta \bm \Xi^{\top} \Delta \bm \Xi$ to ensure the othornormality of $\bm \Xi + \Delta \bm\Xi$.  
%
The perturbation of $\sigma$ is conducted by drawing $d$ random samples (one for each direction) from a uniform distribution $\mathcal{U}(r-\beta, r+\beta)$ with $\beta \ll r$. The random perturbation can be triggered by various types of indicators, e.g., the magnitude of the DGS gradient, the number of iterations completed since last perturbation. The DGS-ES method with the standard gradient descent is summarized in Algorithm 1.

\paragraph{Asymptotic consistency.}
The DGS gradient in Eq.~\eqref{dev_smooth_func} is not an estimator for  $\nabla F_{\sigma}(\bm x)$ or $\nabla F(\bm x)$. In fact, it is designed to be used in the nonlocal setting for non-convex optimization. In nonlocal modeling, a common practice is to study the asymptotic consistency between local and nonlocal gradients (see \cite{doi:10.1137/19M1296720}). A more direct question is that ``{\em Does the DGS gradient estimator converge to the local gradient as $\sigma$ approaches to zero?}'' This question can be answered easily for $F(\bm x) \in \mathcal{C}^{1,1}(\mathbb{R}^d)$. In this case, there exists $L> 0 $ such that
 $
 \|\nabla F(\bm x + \bm \xi) - \nabla F(\bm x)\| \le L\|\bm \xi\|, \, 
 \forall \bm x ,\, \bm \xi\in \mathbb{R}^d
 $ ($\|\cdot\|$ denotes the $L^2$ norm in this work). 
 Then, the difference between $\widetilde{\nabla}^M_{ \sigma, \bm \Xi}[F]$ and $\nabla F$ can be bounded by 
 \begin{align*}
& \left\| \widetilde{\nabla}^M_{ \sigma, \bm \Xi}[F] -  {\nabla}F \right\|^2  
 \le \frac{2 C^2 {\pi} d (M\,!)^2}{4^M((2M)\,!)^2}  \sigma^{4M-2} + 32 d L^2  \sigma^2, \notag
\end{align*}
where the first term on the right hand side comes from the GH quadrature and the second term measures the difference between $\nabla F$ and $\nabla_{\sigma, \bm \Xi}[F]$. It is easy to see the asymptotic consistency, i.e., 
\[
\lim_{\sigma \rightarrow 0} \big| \nabla F(\bm x) - \widetilde{\nabla}_{\sigma, \bm \Xi}^M[F](\bm x) \big| = 0
\]%
\begin{wraptable}{r}{0.52\textwidth}
\vspace{-0.2cm}
\footnotesize
\begin{tabular}{p{0.48\textwidth}}\\
\toprule
  {\;\;\bf Algorithm 1: The DGS-ES algorithm} \\\midrule
\begin{algorithmic}[1]
\setstretch{1.12}
  \STATE{\bf Hyper-parameters}:$M$: \# GH quadrature points; $\lambda_t$: learning rate; $\alpha$: the scaling factor for the rotation $\Delta \bm \Xi$; $r, \beta$: the mean and radius for sampling $\sigma$; $\gamma$: the tolerance for triggering random perturbation.\\
  \STATE{\bfseries Input:} The initial state $\bm x_0$
  \STATE{\bfseries Output:} The final state $\bm x_T$
  \STATE Set $\bm \Xi = \mathbf{I}_d$, and $\sigma_i = r$ for $i = 1, \ldots, d$
  \FOR{$t=0, \ldots T-1$}
  \STATE Evaluate $\{G(\sqrt{2}\sigma_i v_m \, | \, \bm x_t, \bm \xi_i)\}^{i=1, \ldots, d}_{m =1, \ldots, M}$
  \FOR{$i = 1, \ldots, d$}
  \STATE Compute $\widetilde{\mathscr{D}}^M[G_{\sigma_i}(0\,|\, \bm x_t, \bm \xi_i)]$ in Eq.~(\ref{e8})
  \ENDFOR
  \STATE Assemble $\widetilde{\nabla}^M_{\sigma, \bm \Xi}[F](\bm x_t)$ in Eq.~(\ref{e5})
  \STATE Set $\bm x_{t+1} = \bm x_t - \lambda_t \widetilde{\nabla}^M_{\bm \sigma, \bm \Xi}[F](\bm x_t) $
  \IF{$\|\widetilde{\nabla}^M_{\sigma, \bm \Xi}[F](\bm x_t)\|_2 < \gamma$}
  \STATE Generate $\Delta \bm \Xi$ and update $\bm \Xi = \mathbf{I}_d + \Delta \bm \Xi$
  \STATE Generate $\sigma_i$ from $\mathcal{U}(r-\beta, r+\beta)$
  \ENDIF
  \ENDFOR
\end{algorithmic}\\\bottomrule
\end{tabular}
\vspace{-1.7cm}
\end{wraptable}
for $M>2$ regardless of the choice of $\bm \Xi$. Consequently, to achieve $\| \widetilde{\nabla}^M_{ \sigma, \bm \Xi}[F] -  {\nabla}F \| \le \varepsilon$ for a fixed $\varepsilon >0$, we need $\sigma \le {\varepsilon}/(4L\sqrt{d})$
and $M \ge \log({2d}/{\varepsilon^2})$, which means the total number of function evaluations should be bigger than $d\log({2d}/{\varepsilon^2})$.
%
We remark that to acquire the same consistency, the MC estimator for $\nabla F_\sigma$ in Eq.~\eqref{e40} requires $\sigma\le {\varepsilon}/(Ld)$ and the number of function evaluations to be $O({d \|{\nabla}F(\bm x) \|^2 }/{\varepsilon^2})$ (see \cite{BCCS19}). Thus, the DGS estimator requires fewer function evaluations to find search direction even in the local setting. Moreover, it seems not easy to relax this dependency of the MC estimator on $d$ and $\varepsilon$ with variance reduction techniques, e.g., \cite{CRSTW18,10.5555/3326943.3326962}.

\section{Experiments}\label{sec:ex}
We present the experimental results using three sets of problems. All experiments were implemented in Python 3.6 and conducted on a set of cloud servers with Intel Xeon E5 CPUs. We compare the DGS-ES method with the following (a) {\bf ES-Bpop}: the standard OpenAI evolution strategy in \cite{SHCS17} with a big population (i.e., using the same number of samples as DGS-ES), (b) {\bf ASEBO}\footnote{We do not report comparison with some recent work on ES methods, e.g., \cite{10.1145/2908812.2908863,8410043}, because they underperform ASEBO as shown in \cite{Choromanski_ES-Active} and code available at \url{https://github.com/jparkerholder/ASEBO}.}: Adaptive ES-Active Subspaces for Blackbox Optimization \cite{Choromanski_ES-Active} with a population of size $4+3\log(d)$, (c) {\bf IPop-CMA}: the restart covariance matrix adaptation evolution strategy with increased population size \cite{1554902}, (d) {\bf Nesterov}: the random search method in \cite{NesterovSpokoiny15}, and (e) {\bf FD}: the classical central difference scheme. The information of the codes used for the baselines is provided in Appendix.



\subsection{Tests on benchmark functions for global optimization}\label{sec:ex_2}
We test the DGS-ES performance on six 2000D benchmark functions \cite{10.1145/1830761.1830794,Jamil2013ALS} for global optimization, i.e., $F_1(\bm x)$: Sphere, $F_2(\bm x)$: Sharp Ridge, $F_3(\bm x)$: Ackley, $F_4(\bm x)$: Rastrigin, $F_5(\bm x)$: Schaffer, and $F_6(\bm x)$: Schwefel. Their definitions and properties are described in Appendix. We performed grid search to tune the hyper-parameters for all the algorithms to ensure a fair comparison. Description of the hyper-parameter tuning can also be found in Appendix. 

The results are shown in Figure \ref{fig2} and Table \ref{sample-table1}. DGS-ES has the best performance overall. In particular, DGS-ES demonstrates significantly superior performance in optimizing the highly non-convex functions $F_3$, $F_4$ and $F_5$. We explain the reasons from the two advantages of DGS-ES.  

{\bf \em Nonlocal exploration}. We use the averaged cosine distance in Table \ref{sample-table1} to compare the nonlocal exploration performance of each method in minimizing the six benchmark functions,
\begin{equation}\label{cos_dist}
{\rm Cos\_Dist} = \frac{1}{T}\sum_{t=1}^T \left(1 - \dfrac{\langle \bm x_t-\bm x_{t-1}, \bm x^*-\bm x_{t-1}\rangle}{\|\bm x_t-\bm x_{t-1}\| \|\bm x^*-\bm x_{t-1}\|}\right),
\end{equation}
where $\bm x^*$ is the global minimum and $T$ is the number of iterations of an optimization path\footnote{The Cos\_Dist is generated along different paths for different methods.}. The Cos\_Dist of the test cases in Figure \ref{fig2} are shown in Table \ref{sample-table1}.
We have the following findings: (1) the DGS-ES provides smallest Cos\_Dist in most cases, which demonstrates that the DGS gradient is very close to the direction pointing to the global minimum (even for functions with many local minima, e.g., $F_3$, $F_4$, $F_5$). (2) FD achieves similar performance with DGS-ES for $F_1$ and $F_2$ where the local gradients also point to the global minimum, but FD is trapped in local minima for non-convex $F_3$, $F_4$, $F_5$. (3) As Nesterov randomly selects the search direction, it is reasonable that most search directions are perpendicular to $\bm x^*-\bm x_{t-1}$. 
(4) ES-Bpop, ASEBO and IPop-CMA have too few random samples to capture the optimal direction $\bm x^*-\bm x_{t-1}$.

\begin{figure}
     \centering
  \includegraphics[scale = 0.33]{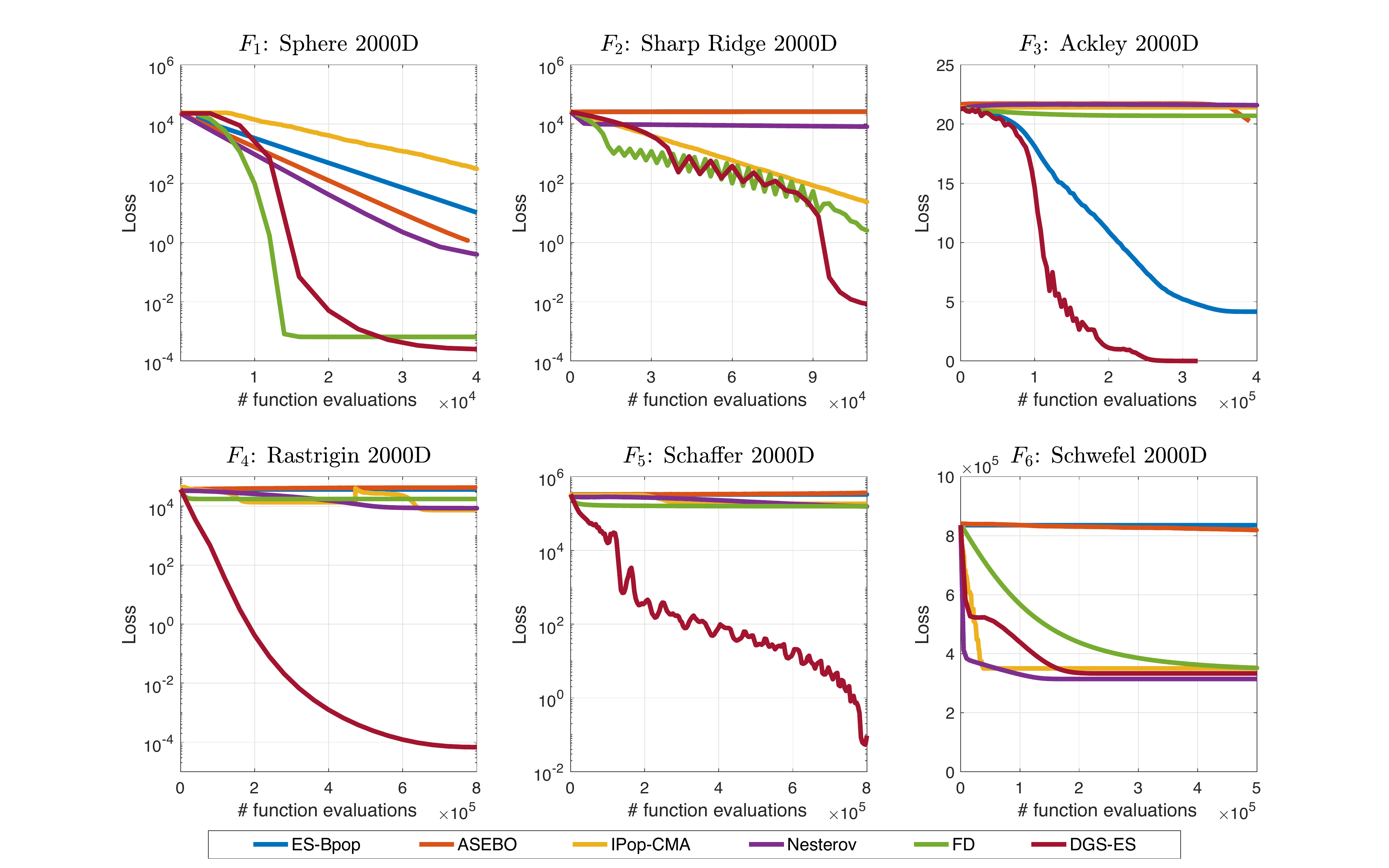}
    \caption{Comparison of the loss decay w.r.t.~\# function evaluations for the 6 benchmark functions in 2000-dimensional spaces. Each curve was generated by averaging 20 independent trials with random initial states. The global minimum is $F(\bm x)=0$ for all the six functions. DGS-ES has the best performance overall, especially for the highly non-convex functions $F_3$, $F_4$, $F_5$. All the methods fail to find the global minimum of $F_6$ which has no global structure to exploit.} 
    \label{fig2}
\end{figure}
\begin{table}[h!]
\vspace{-0.0cm}
  \centering
  \small
  \begin{tabular}{p{1.6cm}p{1.1cm}p{1.5cm}p{1.1cm}p{1.5cm}p{1.1cm}p{1.5cm}}
    \toprule
    &  \multicolumn{2}{c}{{{$F_1$: Sphere}}} & \multicolumn{2}{c}{{{$F_2$: Sharp Ridge}}} & 
    \multicolumn{2}{c}{{{$F_3$: Ackley}}}\\
    \cmidrule(r){2-3} \cmidrule(r){4-5} \cmidrule(r){6-7}
      &  {Cos\_Dist} & {Grad\_Norm} &  {Cos\_Dist} & {Grad\_Norm}&  {Cos\_Dist} & {Grad\_Norm}\\
    \midrule
    DGS-ES       & {\bf 1.86e-9} & 8.09e+1 & 1.48e-1 & 4.11e+1 & {\bf 7.71e-2} & {\bf 5.49e-2}\\
    ES-Bpop      & 2.91e-1 & 1.08e+2 & 3.23e-1 & 3.88e+2 & 7.03e-1 & 1.37e-1\\
    ASEBO        & 9.25e-1 & 8.08e+2 & 9.25e-1 & 2.56e+2 & 1.00e+0 & 8.13e-1\\
    Nesterov     & 9.82e-1 & 2.75e+3 & 9.82e-1 & 2.70e+3 & 9.99e-1 & 3.28e+0\\
    FD           & 1.81e-1 & {\bf 8.04e+1} & {\bf 9.64e-2} & {\bf 2.60e+1} & 9.82e-1 & 7.38e-2\\
    IPop-CMA       & 1.07e+0 & N/A & 1.07e+0 & N/A & 9.99e-1 & N/A\\
    \toprule
    &  \multicolumn{2}{c}{{{$F_4$: Rastrigin}}} & \multicolumn{2}{c}{{{$F_5$: Schaffer}}} & 
    \multicolumn{2}{c}{{{$F_6$: Schwefel}}}\\
    \cmidrule(r){2-3} \cmidrule(r){4-5} \cmidrule(r){6-7}
      &  {Cos\_Dist} & {Grad\_Norm} &  {Cos\_Dist} & {Grad\_Norm}&  {Cos\_Dist} & {Grad\_Norm}\\
    \midrule
    DGS-ES       & {\bf 3.01e-5} & {\bf 5.76e+1} & {\bf 4.85e-1} & {\bf 1.50e+1} & 1.04e+0 & 8.45e+1\\
    ES-Bpop      & 9.27e-1 & 6.46e+2 & 9.19e-1 & 5.06e+2 & 1.00e+0 & 1.23e+3\\
    ASEBO        & 9.93e-1 & 6.27e+3 & 9.98e-1 & 3.88e+4 & {\bf 9.94e-1} & 3.05e+4\\
    Nesterov     & 9.99e-1 & 4.31e+4 & 9.95e-1 & 5.50e+4 & 9.99e-1 & 1.79e+3\\
    FD           & 1.02e+0 & 7.25e+1 & 9.78e-1 & 2.28e+2 & 1.08e+0 & {\bf 6.33e+1}\\
    IPop-CMA     & 1.00e+0 & N/A & 1.01e+0 & N/A & 1.03e+0 & N/A\\
    \bottomrule
    \end{tabular}
      \vspace{0.1cm}
    \caption{The average cosine distance in Eq.~\eqref{cos_dist} and the standard deviation of gradient's $L^2$ norm in Eq.~\eqref{grad_norm} for all the test cases shown in Figure \ref{fig2}. The cosine distance is in the range $[0,2]$. The smaller the Cos\_Dist and the Grad\_Norm, the better the performance of a method.}
  \label{sample-table1}
  \vspace{-0.6cm}
\end{table}
{\bf \em Variation of gradient estimators}. We use the standard deviation of the gradient estimators' $L^2$ norm in Table \ref{sample-table1}, to compare the variance of the estimators,
%
\begin{equation}\label{grad_norm}
    \text{Grad\_Norm} = \left(\frac{1}{T} \sum_{t=1}^T (\|\nabla_t\| - \mu)^2\right)^{1/2}\; \text{ and }\; \mu = \frac{1}{T} \sum_{t=1}^T \|\nabla_t\|,
\end{equation}
where $T$ is the number of iterations and $\nabla_t$ denotes different estimators for different methods\footnote{Since CMA does not have gradient or directional derivative, we don't compare IPop-CMA in this setting.}.
We have the following findings: (1) DGS-ES provides the smallest Grad\_Norm for most test cases due to the good smoothing effect with a large $\sigma$ and the high accuracy of the GH quadrature rule. (2) FD achieves similar
Grad\_Norm as DGS-ES for $F_3$ and $F_4$ because it is trapped in local minima. (3) ES-Bpop and ASEBO cannot well control the Grad\_Norm because of the high variance of the MC-based gradient estimator. (4) The Grad\_Norm of Nesterov is biggest for $F_1,F_2,F_3,F_4,F_5$ because the norm of the gradients of these functions are highly fluctuating and the Nesterov method does not provide as good smoothing as DGS-ES to reduce the fluctuation.

\subsection{Constrained topology optimization for architecture design}\label{sec:ex_3}
We demonstrate the portability of the DGS gradient to constrained optimization using a real-world topology optimization (TO) problem. TO has many applications in engineering \cite{bendsoe2013topology, aage2017giga} and recently attracted attentions in machine learning \cite{hoyer2019neural,yu2019deep,oh2019deep,wu2015system,liu2018narrow,li2020hybrid}. We use DGS-based TO to design a 2D vertical cross section of a bridge from random initial guesses (see Figure \ref{fig4}a). 
\begin{wrapfigure}{r}{0.36\textwidth}
 \vspace{-0.3cm}
\hspace{-0.15cm}\includegraphics[scale = 0.3]{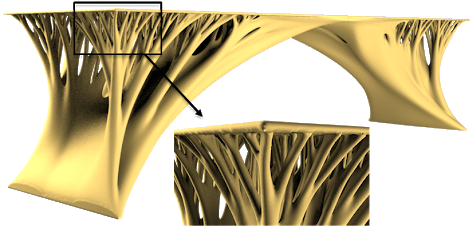}
\caption{Good conceptual design.}\label{fig3}
  \vspace{-0.3cm}
\end{wrapfigure}
The design domain is meshed by $120\times40$ elements, each of which is a design variable ranging from 0 (void) to 1 (solid). By assuming the bridge is symmetric, the total number of independent design variables is $2400~(60 \times 40)$ which is the dimension of the optimization problem. The constraints include (i) $20$\% volume constraint, i.e., the volume of solid materials (black pixels) in Figure \ref{fig4} cannot exceed $480~(2400\times 0.2)$, (ii) unit uniform load on the top and one fixed supports from the bottom. The goal is to optimize the material layout to achieve maximum load-carry capability of the bridge. A conceptually good design is shown in Figure \ref{fig3}.

\begin{figure}[h!]
\centering
\includegraphics[scale = 0.65]{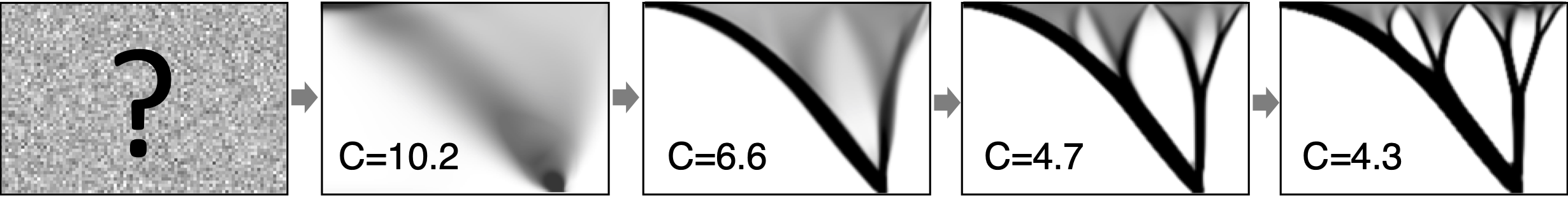}
\caption{\footnotesize Illustration of DGS-based TO design process from random initial guess. The topology of the bridge architecture tends to be more and more clear as the loss function value $C(\times10^5)$ decreases.}\label{fig4}
\end{figure}
%
\begin{figure}[h!]
\vspace{-0.4cm}
\parbox{\textwidth}{
\centering
\begin{minipage}{.6\textwidth}
\centering
  \begin{center}
\includegraphics[width=0.95\textwidth]{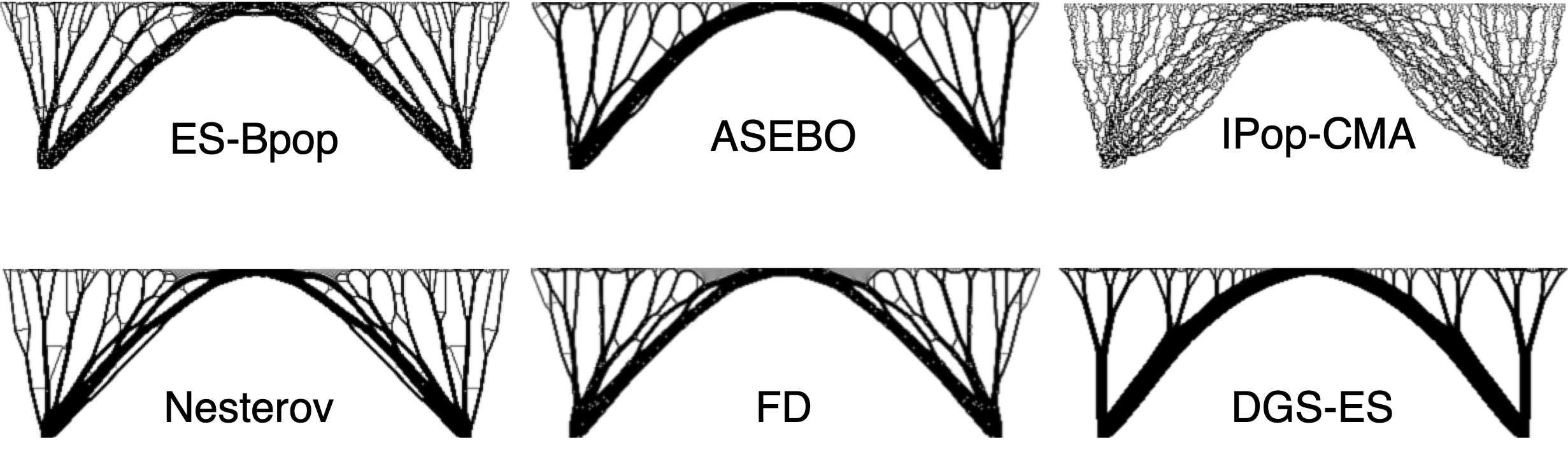}
  \end{center}
\caption{\small Comparison of final typologies. The DGS-based design shows a strong hierarchical tree feature that matches the conceptual design in Figure \ref{fig3}. IPop-CMA tends to a blurry topology. The other algorithms show many local/minor features that have negative impacts on load-carry capability and bridge construction.}\label{fig5} 
\end{minipage}
\begin{minipage}{.39\textwidth}
\hspace{-0.1cm}\includegraphics[scale=0.42]{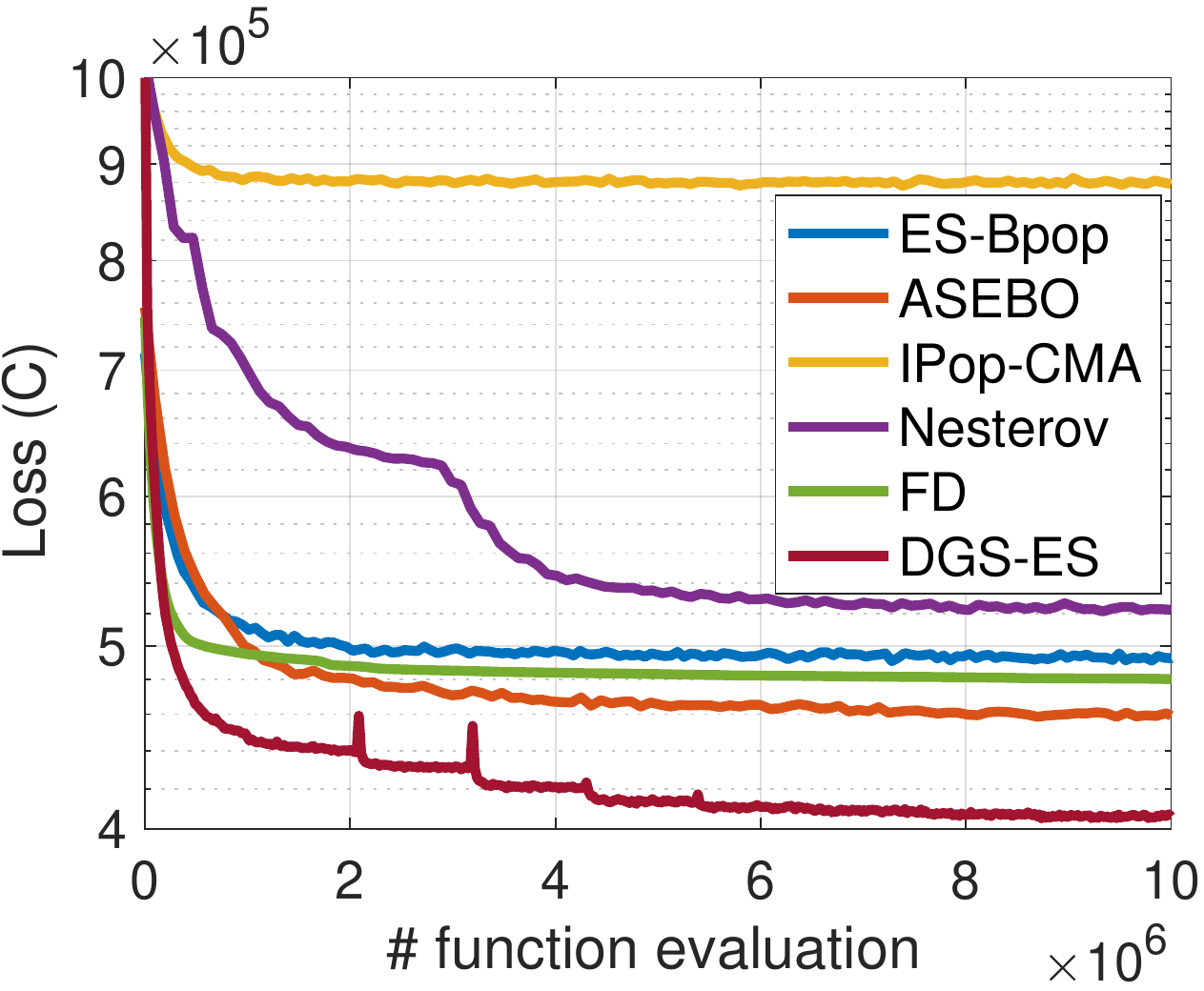}
\vspace{-0.1cm}
\caption{\small Loss decay for DGS-based TO.}\label{fig6}
\end{minipage}
}
\end{figure}
\vspace{-0.1cm}
%
The challenges in TO include highly non-convex and multi-modal loss functions and rigid constraints. Extensive research efforts have been made on developing exclusive constrained optimization algorithms for TO. The state-of-the-art is Method of Moving Asymptotes  (MMA)\cite{svanberg1987method}, which is a gradient-based method. 
%
%
However, MMA is limited to seek optima using local gradients, either via adjoint method or FD. 
Here, we address this issue by inserting the DGS gradient into the MMA framework, and exploit the nonlocal exploration ability of the DGS gradient to find a better design. The hyper-parameters for the DGS gradient will be given in Appendix. 
Figure \ref{fig4} shows the iterative optimization procedure using the DGS-based MMA optimizer.


%

Figure \ref{fig6} summarizes the results\footnote{All the baselines except for IPop-CMA can be inserted into MMA, where the source code can be found at \url{https://github.com/arjendeetman/TopOpt-MMA-Python}.}. We ran each algorithm for 5 times with random initial guesses and plot the mean loss decay. The DGS gradient leads to faster convergence and better final design than the baselines. FD converges fast initially but is quickly trapped into a local minimum. ASEBO and ES-Bpop perform similar to FD. Nesterov may perform well eventually but converges slowly. The IPop-CMA has the worst performance because the simple Lagrangian penalty is insufficient to enforce the constraints. 
%
%
The performances are also demonstrated by the final topology in Figure \ref{fig5}.

\subsection{Inference of hydraulic conductivity field in subsurface environments}\label{sec:ex_4}
\begin{wrapfigure}{r}{0.6\textwidth}
\vspace{-0.6cm}
   \includegraphics[scale = 0.3]{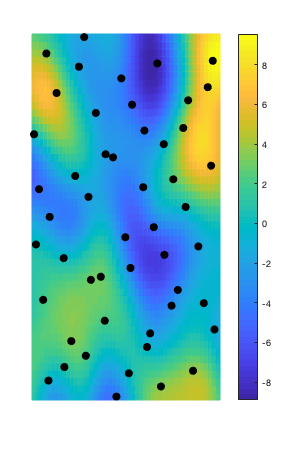}
 \hspace{-0.4cm} \includegraphics[scale = 0.3]{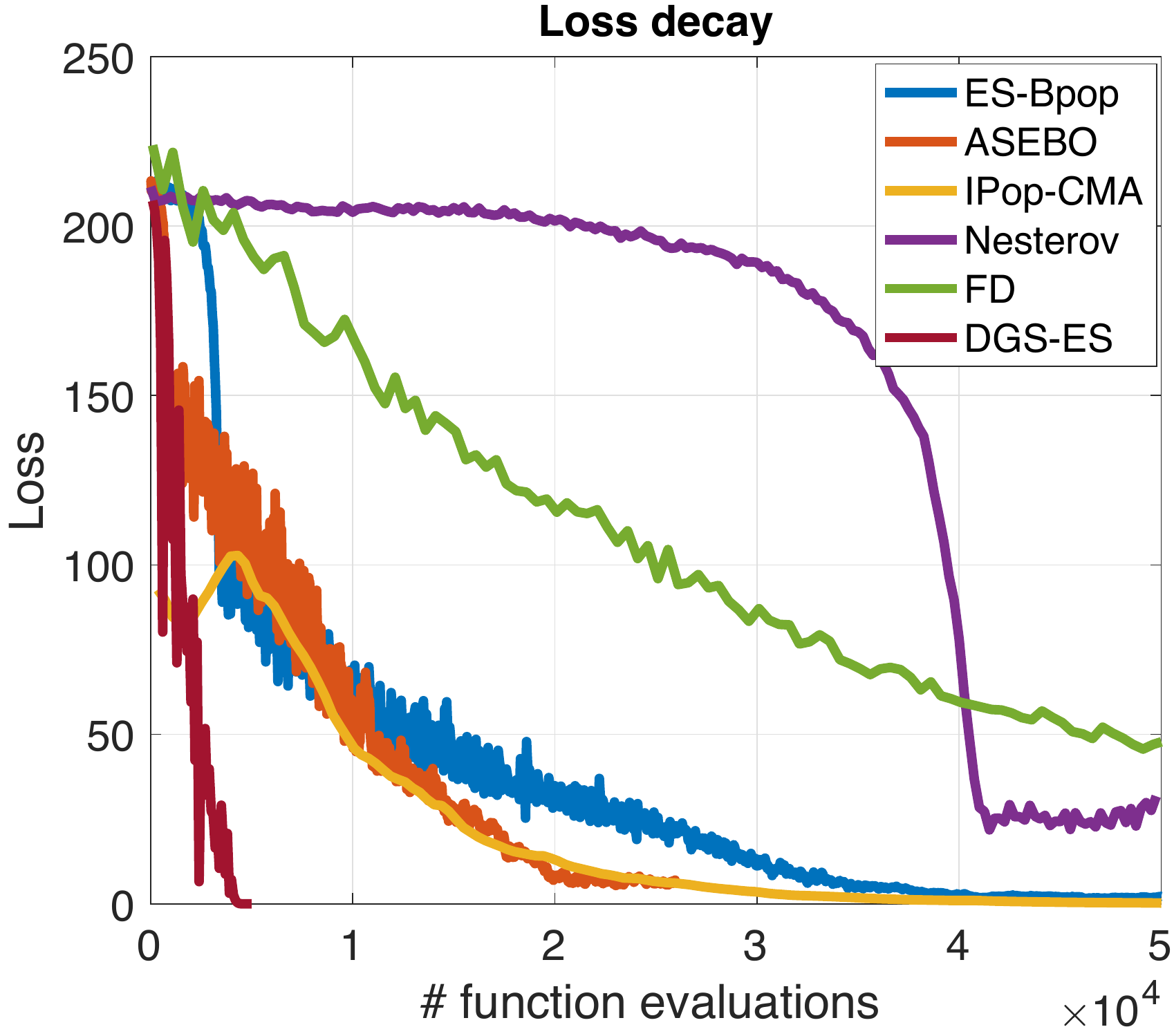}
  \vspace{-0.1cm}
    \caption{(Left): the target hydraulic conductivity field and the 50 locations (black dots) for collecting hydraulic head data. (Right): comparison of the loss decay w.r.t.~\# function evaluations for predicting the hydraulic conductivity field using hydraulic head data.}
    \label{fig7}
    \vspace{-0.3cm}
\end{wrapfigure}
We demonstrate the superior performance of DGS-ES in solving an inference problem in groundwater modeling. Hydraulic conductivity measuring the ease of liquid flow through porous media is an important parameter in predicting contaminant transport in groundwater. However, hydraulic conductivity is very difficult to measure and typically inferred from hydraulic heads (easier to measure). In this work, we use a fully connected neural network (FNN) to approximate a 2D hydraulic conductivity field (Figure \ref{fig7} (left)). The FNN has one hidden layer with 64 neurons. The input is the 2D spatial coordinates and the output is hydraulic conductivity values. $tanh(\cdot)$ is used as the activation. The training data are hydraulic head samples randomly selected at 50 locations. To map the output of the FNN to the training data space, we need to run a blackbox groundwater simulator MODFLOW \cite{modflow} which solves a second-order parabolic partial differential equation. The loss function is defined as the mean squared error between the predicted hydraulic heads and the training data\footnote{The training data is generated by running MODFLOW with the true hydraulic conductivity field.}. The hyper-parameters for all the methods and the parameter values for MODFLOW are given in Appendix. The results in Figure \ref{fig7} (right) clearly demonstrate the much faster convergence of our DGD-ES method compared to other baselines.

\section{Conclusion and discussion}
High-dimensional black-box optimization is an important topic in several machine learning areas, such as reinforcement learning (RL), variational inference and adversarial attacks. We developed the DGS-ES algorithm that takes a novel nonlocal gradient operator with directional Gaussian smoothing to alleviate several challenges in global optimization. Experiments demonstrated that the DGS-ES outperforms several baseline algorithms on both benchmark functions and two scientific problems.

{\bf Limitations}. We realized that there are several limitations with the current version of the DGS-ES algorithm, including (1) \emph{Naive random perturbation strategy.} As the DGS estimator is deterministic, a more effective random exploration strategy is critical to the robustness of the algorithm. (2) \emph{Hyper-parameter tuning.} The most sensitive hyper-parameter is the smoothing radius $\sigma$. If $\sigma$ is too small, the loss function will be insufficiently smoothed, such that the optimizer may be trapped in a local minimum. In contrast, if $\sigma$ is too big, the loss function is overly smoothed, such that the convergence will become much slower. How to adaptively adjust the smoothing radius is still an open question. 
(3) {\em Sub-optimal solution for loss functions without global structures.} Figure \ref{fig2} shows that the DGS-ES cannot find the global minimum of the Schwefel function that does not have a global structure. This could happen in real-world applications. For example, even though our method outperforms the baselines in solving the TO problem, we cannot verify that the design obtained by our method is globally optimal.   

{\bf Future work}. We plan to implement a distributed DGS-ES version to accelerate the time to solution for computationally expensive black-box training problems and demonstrate its strong scalability and the dimension independence property on distributed deep learning frameworks, such as Ray \cite{moritz2018ray}. We also plan on extending the scalable DGS-ES to RL research and help reduce the RL training cost from days to hours or even minutes. That would be a major improvement for the whole RL community.

\section*{Appendix}
\appendix

\renewcommand{\thetable}{\Alph{section}.\arabic{table}}
\renewcommand{\thefigure}{\Alph{section}.\arabic{figure}}

\section{Additional information on the test cases in \S 4.1}

\subsection{Definitions and features of the benchmark functions}\label{app_benchmark}

\begin{itemize}[leftmargin=12pt]\itemsep0.25cm
\item The {\bf Sphere} function $F_1(\bm x)$ is defined by 
%
\begin{equation*}
F_1(\bm{x}) = \sum_{i=1}^d x_i^2,    
\end{equation*}
where $d$ is the dimension and $\bm{x} \in [-5.12, 5.12]^d$ is the input domain. The global minimum is $f(\bm{x}^*) = 0$ at $\bm{x}^* = (0,...,0)$. It represents {\em convex and isotropic} landscapes.

\item The {\bf Sharp Ridge} function $F_2(\bm x)$ is defined by
\begin{equation*}
F_2(\bm x) = x_1^2 + 100 \sqrt{\sum_{i=2}^d x_i^2},
\end{equation*}
where $d$ is the dimension and $\bm{x} \in [-10, 10]^d$ is the input domain. The global minimum is $f(\bm{x}^*) = 0$ at $\bm{x}^* = (0,...,0)$. This represents {\em convex and anisotropic} landscapes. There is a sharp ridge defined along $x_2^2 + \cdots + x_d^2 = 0$ that must be followed to reach the global minimum, which creates difficulties for optimizations algorithms.

\item The {\bf Ackley} function $F_3(\bm x)$ is defined by
\begin{equation*}
F_3(\bm{x}) = -a\exp\left(-b\sqrt{\frac{1}{d}\sum_{i=1}^d x_i^2} \right)-\exp\left( \frac{1}{d}\sum_{i=1}^d \cos(c x_i)\right) + a + \exp(1),
\end{equation*}
where $d$ is the dimension and $a=20, b= 0.2, c=2\pi$ are used in our experiments. The input domain $\bm{x} \in [-32.768, 32.768]$. The global minimum is $f(\bm{x}^*) = 0$, at $\bm{x}^* = (0,...,0)$. The Ackley function represents {\em non-convex} landscapes with {\em nearly flat outer region}.  The function poses a risk for optimization algorithms, particularly hill-climbing algorithms, to be trapped in one of its many local minima.

\item The {\bf Rastrigin} function $F_4(\bm x)$ is defined by
\begin{equation}
    F_4(\bm{x}) = 10d + \sum_{i=1}^d [x_i^2 - 10\cos(2\pi x_i)],
\end{equation}
where $d$ is the dimension and $\bm{x} \in [-5.12, 5.12]^d$ is the input domain. The global minimum is $f(\bm{x}^*) = 0$ at $\bm{x}^* = (0,...,0)$. This function represents {\em multimodal and separable} landscapes.

\item The {\bf Schaffer} function $F_5(\bm x)$ is defined by 
\begin{equation*}
F_5(\bm x) = \frac{1}{d-1} \left(\sum_{i=1}^{d-1}\left(\sqrt{s_i} + \sqrt{s_i} \sin^2(50 s_i^{\frac{1}{5}})\right)\right)^2 \quad \text{with}\quad s_i = \sqrt{x_i^2 + x_{i+1}^2},
\end{equation*}
where $\bm{x} \in [-100, 100]^d$ is the input domain. The global minimum is $f(\bm{x}^*) = 0$, at $\bm{x} = (0,\cdots,0)$. This function represents {\em multimodal and non-separable} landscapes.

\item The {\bf Schwefel} function $F_6(\bm x)$ is defined by
\begin{equation*}
F_6(\bm x) = 418.9829d - \sum_{i=1}^d x_i \sin(\sqrt{|x_i|}),
\end{equation*}
where $\bm{x} \in [-500, 500]^d$ is the input domain. The global minimum is $f(\bm{x}^*) = 0$, at $\bm{x} = (420.9687,\cdots,420.9687)$. This function represents {\em multimodal} landscapes with {\em no global structure}.
\end{itemize}

\subsection{Experimental details}
As we mentioned in the main paper, we perform grid search to tune the hyper-parameters for all algorithms to ensure a fair comparison. Now we describe the hyper-parameter tuning process for each method for the six benchmark functions.

\subsubsection{The DGS-ES method}
There are six hyper-parameters given in Algorithm 1, i.e., $M$: the number of GH quadrature points; $\lambda_t$: learning rate; $\alpha$: the scaling factor for the rotation $\Delta \bm \Xi$; $r, \beta$: the mean and variation for sampling $\sigma$; $\gamma$: the tolerance for triggering random perturbation. Our preliminary study shows that the random perturbation strategy proposed in \S \ref{sec:ada_DGS-ES} does not help improve the performance of the DGS-ES method for the six functions, so that we do not turn on the random perturbation in this example. The other parameters are tuned as follows. 

We tune $M$ by the grid $\{3,5,7,9,11,13,15,17,19,21\}$. The reason why choosing odd numbers is that the center quadrature point can be reused for all the directions to saving computational cost. We tune the learning rate $\lambda_t$ by fitting it to a polynomial decay schedule
\begin{equation}\label{lr_schedule}
\lambda_t = (\lambda_0-\lambda_T) \left(1-\frac{t}{T}\right)^\tau + \lambda_T,
\end{equation}
where $T$ is the maximum number of iterations, $\lambda_0$ is the initial learning rate, $\lambda_T$ is the final learning rate, and $\tau$ is the decay power. We tune the learning rate schedule by finding appropriate parameters in Eq.~\eqref{lr_schedule}. We tune the smoothing radius\footnote{$r=\sigma$ when having no perturbation} $\sigma$ is also tuned by fitting it to a polynomial decay schedule
\begin{equation}\label{sig_schedule}
\sigma_t = (\sigma_0-\sigma_T) \left(1-\frac{t}{T}\right)^\nu + \sigma_T,
\end{equation}
where $T$ is the maximum number of iterations, $\sigma_0$ is the initial mean radius, $\sigma_T$ is the final mean radius, and $\nu$ is the decay power. We tune the smoothing radius schedule by finding appropriate parameters in Eq.~\eqref{sig_schedule}. We report the tuned hyper-parameters for DGS-ES in Table \ref{sample-table2}.
\begin{table}[h!]
\centering
\caption{The hyper-parameter values for DGS-ES used in \S 4.1}
\label{sample-table2}
\begin{tabular}{lcccccccc}
\toprule
            & $M$  & $\lambda_0$ & $\lambda_T$ & $\tau$ & $r_0$ & $r_T$ & $\nu$ & $T$ \\
\midrule
Sphere 2000D      &  3   & 1.0  & 0.01   &  2.0  & 1.0  & 0.0001 & 2.0 & 10\\
Sharp-Ridge 2000D &  3   & 0.4  & 0.0001 &  3.0  & 0.5  & 0.1    & 0.5 & 30\\
Ackley 2000D      &  3   & 8000.0  & 0.001 &  4.0  & 2.0  & 0.001 & 2.0 & 80\\
Rastrigin 2000D   &  21  & 0.5  & 0.001 &  2.0  & 1.0  & 0.5 & 2.0 & 20\\
Schaffer 2000D   &  3  & 5.0  & 0.001 &  1.0  & 50.0  & 0.001 & 2.0 & 200\\
Schwefel 2000D   &  5  & 10.0  & 1.0 &  1.0  & 5.0  & 1.0 & 2.0 & 125\\
\midrule
Sphere 20D      &  3   & 1.0  & 0.01   &  2.0  & 1.0  & 0.0001 & 2.0 & 10\\
Sharp-Ridge 20D &  3   & 0.4  & 0.00001 &  4.0  & 0.5  & 0.001    & 2.0 & 30\\
Ackley 20D      &  3   & 200.0  & 0.01 &  4.0  & 2.0  & 0.01 & 2.0 & 80\\
Rastrigin 20D   &  21  & 0.5  & 0.001 &  2.0  & 1.0  & 0.5 & 2.0 & 10\\
Schaffer 20D   &  3  & 5.0  & 0.0001 &  3.0  & 10.0  & 0.001 & 2.0 & 200\\
Schwefel 20D   &  5  & 10.0  & 1.0 &  2.0  & 5.0  & 1.0 & 2.0 & 200\\
\bottomrule
\end{tabular}
\end{table}

\subsubsection{The ES-Bpop method}
ES-Bpop refers to the standard OpenAI evolution strategy in \cite{SHCS17} with the a big population, i.e., the same population size as the DGS-ES method. 
The purpose of using a big population is to compare the MC-based estimator for the standard GS gradient and the GH-based estimator for the DGS gradient given the same computational cost. In this setting, the ES-Bpop only has two hyper-parameters, i.e., the learning rate $\lambda_t$ and the smoothing radius $\sigma_t$. We use the same type of schedule as in Eq.~\eqref{lr_schedule} and Eq.~\eqref{sig_schedule} to tune the two hyper-parameters and the tuned hyper-parameters are given in Table \ref{sample-table3}.
\begin{table}[h!]
\centering
\caption{The hyper-parameter values for ES-Bpop used in \S 4.1}
\label{sample-table3}
\begin{tabular}{lccccccc}
\toprule
            &  $\lambda_0$ & $\lambda_T$ & $\tau$ & $r_0$ & $r_T$ & $\nu$ & $T$ \\
\midrule
Sphere 2000D       & 0.1  & 0.01   &  2.0  &  0.001  & 0.0001 & 2.0 & 20\\
Sharp-Ridge 2000D  & 0.001  & 0.0001 &  2.0  & 0.5  & 0.1    & 0.5 & 30\\
Ackley 2000D       & 1000.0  & 0.001 &  4.0  & 2.0  & 0.001 & 2.0 & 100\\
Rastrigin 2000D    & 0.01  & 0.001 &  2.0  & 1.0  & 0.5 & 2.0 & 20\\
Schaffer 2000D   & 0.2  & 0.1 &  2.0  & 50.0  & 0.001 & 2.0 & 200\\
Schwefel 2000D   & 0.0001  & 0.00001 &  2.0  & 1.0  & 0.1 & 2.0 & 125\\
\midrule
Sphere 20D        & 0.5  & 0.01   &  3.0  & 1.0  & 0.01 & 2.0 & 10\\
Sharp-Ridge 20D   & 0.1  & 0.0001 &  2.0  & 1.0  & 0.001    & 2.0 & 30\\
Ackley 20D        & 50.0  & 0.01 &  4.0  & 2.0  & 0.01 & 2.0 & 200\\
Rastrigin 20D     & 0.1  & 0.001 &  2.0  & 1.0  & 0.1 & 2.0 & 15\\
Schaffer 20D     & 0.5  & 0.1 &  2.0  & 10.0  & 1.0 & 2.0 & 200\\
Schwefel 20D     & 0.5  & 0.1 &  2.0  & 1.0  & 0.5 & 2.0 & 200\\
\bottomrule
\end{tabular}
\end{table}

\subsubsection{The ASEBO method}
{ASEBO} refers to Adaptive ES-Active Subspaces for Blackbox Optimization proposed in \cite{Choromanski_ES-Active}. This is the state-of-the-art method in the family of ES. It has been shown that other recent developments on ES, e.g., \cite{10.1145/2908812.2908863,8410043}, underperform ASEBO in optimizing the benchmark functions. We use the code published at \url{https://github.com/jparkerholder/ASEBO} by the authors of the ASEBO method. Since we use ASEBO to represent the state-of-the-art ES method, we set the population size to the standard value $4+3\log(d)$. Due to the use of a small population, the smoothing radius $\sigma$ needs to be small to control the variance of the MC estimator. It turns out that $\sigma=0.1$ works the best for the test cases. The ASEBO code uses $\lambda_t = \alpha \lambda_{t-1}$ to set a learning rate decay schedule. For the six test cases, we tune $\alpha$ by searching the grid $\{0.999,0.99,0.9\}$ and choose $\alpha =0.99$ for all the six functions. Again, due to the small population size, the initial learning rate $\lambda_0$ cannot be set to as large as those for DGS-ES or ES-Bpop because it will overshoot. Thus, we search for the initial learning rate on the grid $\{0.01,0.1,0.5,1.0\}$, and it turns out $\lambda_0=0.1$ provides the best performance for ASEBO for the test functions. 

\subsubsection{The IPop-CMA method}
IPop-CMA refers to the restart covariance matrix adaptation evolution strategy with increased population size proposed in \cite{1554902}. We use the code pycma v3.0.3 available at \url{https://github.com/CMA-ES/pycma}. The main subroutine we use is \texttt{cma.fmin}, in which the hyper-parameters are 
\vspace{-0.2cm}
\begin{itemize}[leftmargin=12pt]
    \item \texttt{restarts=9}: the maximum number of restarts with increasing population size;
    \item \texttt{restart\_from\_best=False}: which point to restart from;
    \item \texttt{incpopsize=2}: multiplier for increasing the population size before each restart;
    \item \texttt{$\sigma_0$}: the initial exploration radius is set to 1/4 of the search domain width.
\end{itemize}

\subsubsection{The Nesterov method}
Nesterov refers to the random search method proposed in \cite{NesterovSpokoiny15}. We use the stochastic oracle 
\[
\bm x_{t+1} = \bm x_{t} - \lambda_t F'(\bm x_t, \bm u_t),
\]
where $\bm u_t$ is a randomly selected direction and $F'(\bm x_t, \bm u_t)$ is the directional derivative along $\bm u_t$. According to the analysis in \cite{NesterovSpokoiny15}, this oracle is more powerful and can be used for non-convex non-smooth functions. As suggested in \cite{NesterovSpokoiny15}, we use forward difference scheme to compute the directional derivative. The only hyper-parameter is the learning rate $\lambda_t$. We use the same polynomial decay model as in DGS-ES and the tuned learning rate schedule is given in Table \ref{sample-table4}
\begin{table}[h!]
\centering
\caption{The hyper-parameter values for Nesterov in \S 4.1}
\label{sample-table4}
\begin{tabular}{lcccccccc}
\toprule
              & $\lambda_0$ & $\lambda_T$ & $\tau$  & $T$ \\
\midrule
Sphere 2000D      & 0.001  & 0.0001   &  2.0 & 1,000,000\\
Sharp-Ridge 2000D  & 0.001  & 0.0001 &  2.0  & 1,000,000\\
Ackley 2000D      & 1.0  & 0.001 &  2.0  & 1,000,000\\
Rastrigin 2000D   & 0.00001  & 0.000001 &  2.0  & 1,000,000\\
Schaffer 2000D  & 0.0001  & 0.0001 &  2.0   & 1,000,000\\
Schwefel 2000D   & 0.005  & 0.0001 &  2.0   & 1,000,000\\
\midrule
Sphere 20D      & 0.01  & 0.001   &  2.0 & 800\\
Sharp-Ridge 20D  & 0.001  & 0.00001 &  2.0   & 2000\\
Ackley 20D       & 0.1  & 0.01 &  2.0   & 8000\\
Rastrigin 20D    & 0.001  & 0.0001 &  2.0   & 6000\\
Schaffer 20D  & 0.005  & 0.001 &  2.0   & 8000\\
Schwefel 20D    & 0.1  & 0.01 &  2.0   & 10000\\
\bottomrule
\end{tabular}
\end{table}

\subsubsection{The FD method}
FD refers to the classical central difference scheme for local gradient estimation. The only hyper-parameter is the learning rate $\lambda_t$. We use the same polynomial decay model as in DGS-ES and the tuned learning rate schedule is given in Table \ref{sample-table5}.
\begin{table}[h!]
\centering
\caption{The hyper-parameter values for FD in \S 4.1}
\label{sample-table5}
\begin{tabular}{lcccccccc}
\toprule
              & $\lambda_0$ & $\lambda_T$ & $\tau$  & $T$ \\
\midrule
Sphere 2000D      & 1.0  & 0.01   &  2.0 & 20\\
Sharp-Ridge 2000D  & 0.4  & 0.0001 &  3.0  & 60\\
Ackley 2000D      & 1.0  & 0.001 &  2.0  & 160\\
Rastrigin 2000D   & 0.001  & 0.0001 &  2.0  & 400\\
Schaffer 2000D  & 0.01  & 0.001 &  2.0   & 600\\
Schwefel 2000D   & 0.1  & 0.01 &  2.0   & 500\\
\midrule
Sphere 20D      & 1.0  & 0.01   &  2.0 & 40\\
Sharp-Ridge 20D  & 0.4  & 0.0001 &  4.0   & 50\\
Ackley 20D       & 1.0  & 0.01 &  4.0   & 400\\
Rastrigin 20D    & 0.01  & 0.001 &  2.0   & 300\\
Schaffer 20D  & 0.001  & 0.0001 &  2.0   & 400\\
Schwefel 20D    & 0.1  & 0.01 &  2.0   & 500\\
\bottomrule
\end{tabular}
\end{table}

\subsection{Additional results and discussion}
Besides high-dimensional problems, it is also natural to compare the performance of the DGS-ES method with the other baselines in solving relatively low-dimensional problems. To this end, we conduct another set of tests for the functions in Appendix A.1 in 20-dimensional spaces. The hyper-parameters used for the 20D tests are also given in Appendix A.2. The main results are given in Figure \ref{fig10}. We have the following observations by comparing the 2000D and 20D results and hyper-parameter values used for DGS-ES and the baselines.
\vspace{-0.1cm}
\begin{itemize}[leftmargin=13pt]
    \item The DGS-ES method can use more aggressive learning rate schedules to accelerate convergence. This is due to the good smoothing effect of the DGS gradient and the high accuracy of the DGS estimator. For example, the Ackley function has a very large and flat outer-region, and the convergence speed depends on how fast an optimizer can go through the flat region and get to the mode containing the global minimum. DGS-ES can use a very large initial learning rate e.g., $\lambda_0=8000$ for the 2000D case, without worry about overshooting, because its cosine distance (shown in Table 1) is small and the variance of the DGS estimator (shown in Table 1) is also small. In comparison, ES-Bpop, Nesterov and FD need to use much smaller learning rates to avoid overshooting; and ES-Bpop performs better than Nesterov and FD due to the use of a relatively large population. Another reason why Nesterov and FD need to use smaller learning rates is that the local derivatives/gradients of some test functions are very fluctuating, e.g., Schaffer and Rastrigin, and both methods do not provide sufficient smoothing effect to reduce the fluctuation.  

    \item The advantage of the DGS-ES in the 20D case is not as significant as in the 2000D case. For example, IPop-CMA outperforms DGS-ES in minimizing the 20D Ackley function. The convergence speed of IPop-CMA type methods depends on how fast samples can be dropped in the mode containing the global minimum. The flat outer-region of 20D Ackley is much smaller than that of the 2000D Ackley, so that it is easier for IPop-CMA to have a sample dropped in the mode containing the global minimum. For the Sharp Ridge function, FD outperforms DGS-ES because the sharp ridge defined by $x_2^2 + \cdots + x_d^2 = 0$ is easier to follow in the 20D space. 
    
    \item The landscapes without any global structures, e.g., the Schwefel function, is still difficult to minimize in 20D cases. 
\end{itemize}
\begin{figure}[h!]
     \centering
  \includegraphics[scale = 0.35]{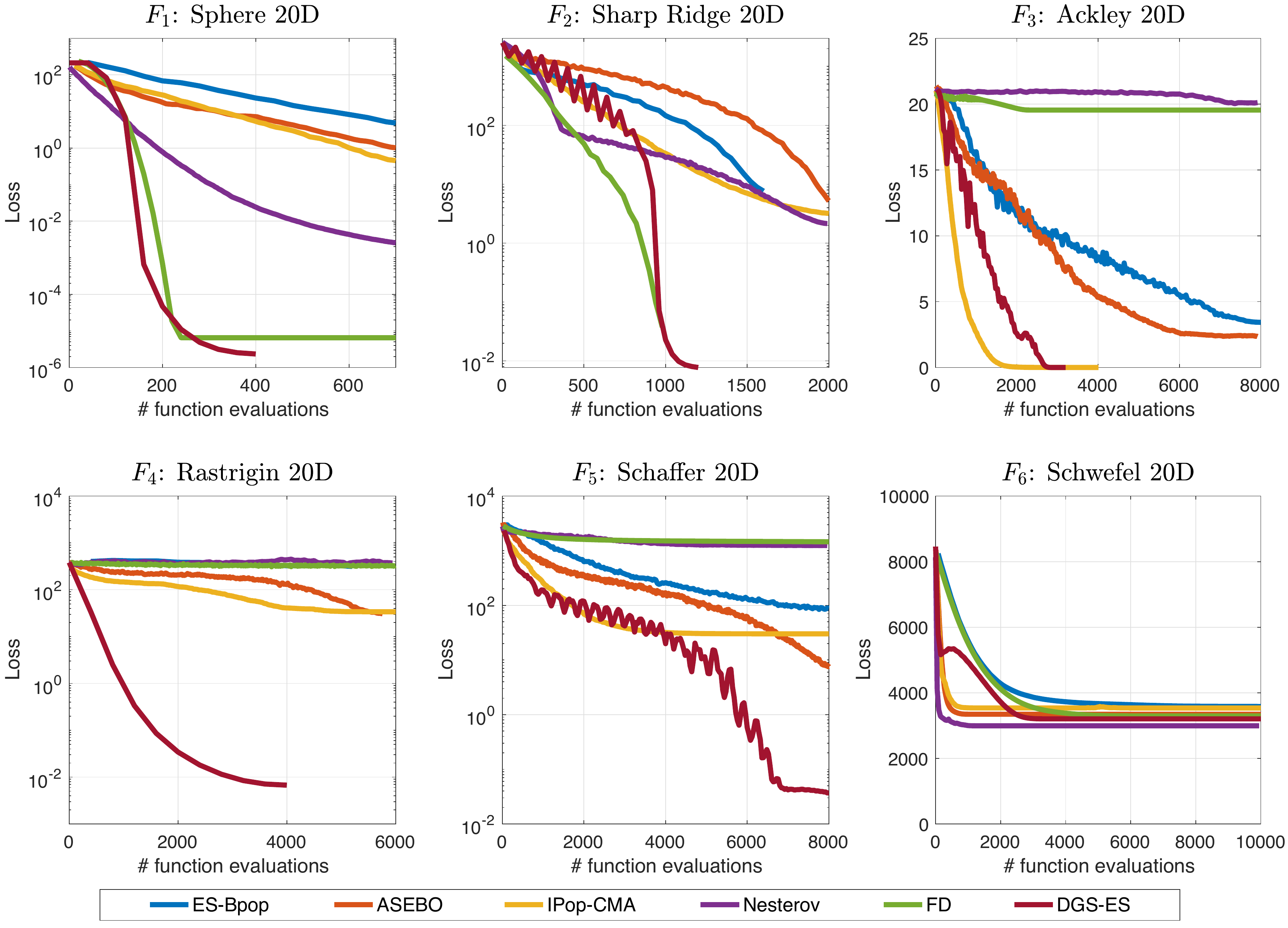}
    \caption{Comparison of the loss decay w.r.t.~\# function evaluations for the 6 benchmark functions in 20-dimensional spaces. Each curve was generated by averaging 20 independent trials with random initialization. The global minimum is $F(\bm x)=0$ for all the six functions.}
    \label{fig10}
\end{figure}

\section{Additional information on the constrained topology optimization in \S 4.2}
\subsection{Topology optimization mathematical formulation}
Here we provide more background information about the topology optimization (TO) problem tested in \S 4.2.
TO is a mathematical method that aims to optimize material layout defined on a design domain $\Omega$ with given boundary conditions, loads and volume constraint, to minimize structural compliance $C$, or equivalently, the least strain energy. In this work, we use the modified Solid Isotropic Material with Penalization (SIMP) approach \cite{sigmund2007morphology} with design-based approach to topology optimization, where each element $e$ is assigned a density $x_e$ that determines its Young's modulus $E_e$:
\begin{equation}
    E_e(x_e) = E_{\min} + x_e^p(E_0-E_{\min}),\quad x_e \in [0,1]
\end{equation}
where $E_0$ is the stiffness of the material, $E_{\min}$ is a very small stiffness assigned to void regions to prevent the stiffness matrix becoming singular. The modified SIMP approach differs from the classical SIMP approach \cite{bendsoe1989optimal}, where elements with zero stiffness are avoided by using a small value. The modified mathematical formulation of the considered TO problem is
\begin{equation}\label{eq: simp2}
\begin{aligned}
\min_{{\bm x}} &: \;\; C(\bm x) = \mathbf{U}^T \mathbf{K} \mathbf{U} = \sum_{e=1}^N E_e(x_e) \mathbf{u}_e^T \mathbf{k}_0 \mathbf{u}_e \\
s.t. &: \;\;{V({\bm x})}/{{V_0}} \le \zeta \\
&: \;\; \mathbf{KU=F} \\
&: \;\; 0 \le \bm x \le 1
\end{aligned}
\end{equation}
where $\bm x$ is the vector of design variables, $C$ is the structural compliance, $\mathbf{K}$ is the global stiffness matrix, $\mathbf{U}$ and $\mathbf{F}$ are the global displacement and force vectors respectively, $\mathbf{u}_e$ and $\mathbf{k}_e$ are the element displacement vector and stiffness matrix respectively, $N$ is the number of elements used to discretize the design domain $\Omega$, $V(\bm{x})$ and $V_0$ are the material volume and design domain volume respectively, $\zeta$ is the prescribed volume fraction, and $p$ is the penalization power coefficient (typically $p=3$). 

\subsection{Topology optimization problem setting in \S 4.2}
\begin{figure}[!h]
\centering
\includegraphics[width=0.4\textwidth]{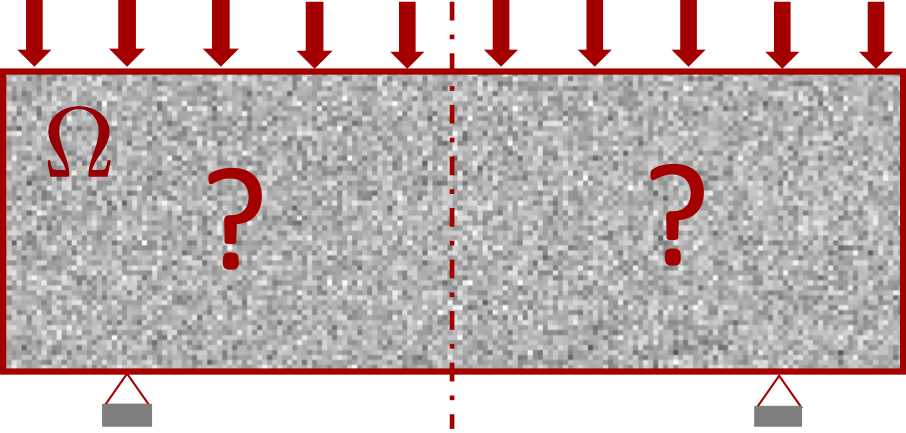}
\caption{Illustration of the structural bridge design}
\label{fig:app_topo1}
\end{figure}

The experimental TO example in \S 4.2 is a typical structural bridge design problem. Our goal is to minimize the structural compliance subject to unit uniform pressure on the top of the design domain $\Omega$, and two fixed supports on the bottom of the design domain $\Omega$, as shown in Fig.\ref{fig:app_topo1}. The volume fraction is $\zeta =0.2$. The design domain $\Omega$ is discretized by 120$\times$40 elements. Using symmetric boundary condition, we reduce the whole design domain into an half to save computational cost. As a result, there are total 2400 (60$\times$40) design variables.

There has been several mathematical challenges in solving TO problem, as shown in Eq.\eqref{eq: simp2}, which can be briefly summarized as follows:
\vspace{-0.1cm}
\begin{itemize}[leftmargin=13pt]
    
    \item {\em Rigid constrained optimization}. Typically, design criteria are specified by user to satisfy multiple constraints including material volume fraction, maximum stress and geometry constraints, so that meaningful structures can be obtained.
    
    \item {\em High dimensional design space}. This is because each element is identified as an independent design variable in TO, shown in Eq.\eqref{eq: simp2}. In most cases, it requires a large number of elements to ensure the accuracy of finite element method (FEM) and to perform a clear final topology. 
    
    \item {\em Highly nonconvex propriety}. TO introduces the SIMP approach to convert the 0-1 integer optimization into continuous optimization but it also changes TO to a difficult nonconvex optimization problem, which has been demonstrated by \cite{stolpe2001trajectories, beck2010sequential, bendsoe2013topology}, such that many local minima exist. 
    
    \item {\em Optimum depends on the initial design}. TO has challenges to pursue global minima \cite{bendsoe2013topology} due to a limited capability of global exploration in sensitivity analysis when using adjoint methods or finite difference methods. The different initial guesses therefore finally lead to different local minima.
 
\end{itemize}

These challenges make a few optimization methods infeasible in solving TO problem well. Bayesian optimization (BO) is good at pursuing global minima but has limitations in high-dimensional problems ($D>1000$). While many improvements have been proposed \cite{kandasamy2015high, li2016high, wang2017batched,mutny2018efficient,rana2017high} to mitigate this challenge, the constraints in BO framework is still a critical issue \cite{gardner2014bayesian}, specifically in practical implementation. The widely used training algorithms including SGD, Adam, RMSprop, etc., are inapplicable to TO problems without handling constraints. 

\subsection{The constrained optimization method used in \S 4.2}
Method of Moving Asymptotes (MMA) \cite{svanberg1987method,svanberg2002class} algorithm is the state-of-the-art optimizer, which has been demonstrated to be versatile and well suited for wide range TO problems. The basic of MMA aims at solving general nonlinear constrained optimization problem: 
\begin{equation}
\begin{aligned}
\min_{\mathbf{x}} &: \quad f_0(\mathbf{x}) + a_0z + \sum_{i=1}^m(c_iy_i + \frac{1}{2}d_i y_i^2)\\
s.t. &: \quad f_i(\mathbf{x})-a_i z -y_i \le 0,\quad i = 1,...,m \\
&: \quad \mathbf{x}\in X, \mathbf{y}\ge0, z\ge 0
\label{eq: mma1}
\end{aligned}
\end{equation}
Here, $X = \left\{x \in \mathbb{R}^n | x_j^{\min} \le x_j \le x_j^{\max}, j=1,...,n \right\}$, where $x_j^{\min}$ and $x_j^{\max}$ are given real numbers which satisfy $x_j^{\min} < x_j^{\max}$ for all $j, f_0, f_1, ..., f_m$ are given, continuously differentiable, real-valued functions on $X$, $a_0, a_i, c_i$ and $d_i$ are given real numbers which satisfy $a_0>0, a_i \ge 0, c_i \ge 0$ and $d_i \ge 0$ and $c_i + d_i >0$ for all $i$ and also $a_ic_i >a_0$ for all $i$ with $a_i>0$.  

MMA is a gradient-based method for solving Eq.\eqref{eq: mma1} using the following steps. In each iteration, given the current point $(\mathbf{x}^{(k)}, \mathbf{y}^{(k)},z^{(k)})$, MMA generates an approximating subproblem, where the functions $f_i(\mathbf{x})$ are replaced by convex functions $\hat{f}_i^{(k)}(\mathbf{x})$. The approximating functions are determined by the {\em gradient information} at the current iteration point and moving asymptotes parameters which are updated in each iteration based on information from previous iteration points. The next iteration point $(\mathbf{x}^{(k+1)}, \mathbf{y}^{(k+1)},z^{(k+1)})$ is obtained by solving the subproblem, which looks as follows:
\begin{equation}
\begin{aligned}
\min_{\mathbf{x}} &: \quad \hat{f}_0^{(k)}(\mathbf{x}) + a_0z + \sum_{i=1}^m(c_iy_i + \frac{1}{2}d_i y_i^2)\\
s.t. &: \quad \hat{f}_i^{(k)}(\mathbf{x})-a_i z -y_i \le 0,\quad i = 1,...,m \\
&: \quad \alpha_j^{(k)} \le x_j \le \beta_j^{(k)}, \quad j=1,...,m \\
&: \quad {y}_i\ge0, i = 1,...,m \\
&: \quad z\ge 0
\label{eq: mma2}
\end{aligned}
\end{equation}
where the approximating functions $\hat{f}_i^{(k)}(\mathbf{x})$ are chosen as 
\begin{equation}
\hat{f}_i^{(k)}(\mathbf{x}) = \sum_{j=1}^n \left( \frac{p_{ij}^{(k)}}{u_j^{(k)}-x_j}+ \frac{q_{ij}^{(k)}}{x_j-l_j^{(k)}} \right) + r_i^{(k)}, \quad i=0,1,...,m,
\label{eq: mma3}
\end{equation}
where 
\begin{equation}
p_{ij}^{(k)} = (u_j^{(k)}-x_j^{(k)})^2 \left( 1.001\left(\frac{\partial f_i}{\partial x_j}(\mathbf{x}^{(k)}) \right)^{+}  + 0.001\left(\frac{\partial f_i}{\partial x_j}(\mathbf{x}^{(k)}) \right)^{-} + \frac{10^{-5}}{x_j^{\max} - x_j^{\min}}\right),
\label{eq: mma4}
\end{equation}
\begin{equation}
q_{ij}^{(k)} = (x_j^{(k)}-l_j^{(k)})^2 \left( 0.001\left(\frac{\partial f_i}{\partial x_j}(\mathbf{x}^{(k)}) \right)^{+}  + 1.001\left(\frac{\partial f_i}{\partial x_j}(\mathbf{x}^{(k)}) \right)^{-} + \frac{10^{-5}}{x_j^{\max} - x_j^{\min}}\right),
\label{eq: mma5}
\end{equation}
\begin{equation}
r_i^{(k)} = \hat{f}_i(\mathbf{x}^{(k)}) - \sum_{j=1}^n \left( \frac{p_{ij}^{(k)}}{u_j^{(k)}-x_j}+ \frac{q_{ij}^{(k)}}{x_j-l_j^{(k)}} \right).
\label{eq: mma6}
\end{equation}
Here, $\left(\frac{\partial f_i}{\partial x_j}(\mathbf{x}^{(k)}) \right)^{+}$ denotes the largest of the two numbers $\frac{\partial f_i}{\partial x_j}(\mathbf{x}^{(k)})$ and 0, while $\left(\frac{\partial f_i}{\partial x_j}(\mathbf{x}^{(k)}) \right)^{-}$ denotes the largest of the two numbers $-\frac{\partial f_i}{\partial x_j}(\mathbf{x}^{(k)})$ and 0. 

The central advantage of MMA in TO is the use of {\em separable} and {\em convex} approximations. The separable property means that the necessary conditions of the subproblems do not couple the primary variables and the latter means that dual methods or primal-dual methods can be employed \cite{bendsoe2013topology}. Combined both two can significantly reduce the computational cost needed to solve the subproblems, particularly for problems with multiple constraints. 


\subsection{Experimental details} 
\subsubsection{Implementation of the DGS-ES method in TO}
We solve the TO design problem shown in Fig.~\ref{fig:app_topo1} by inserting the DGS gradient into the MMA optimizer, as discussed in main paper. Our idea is to exploit the nonlocal exploration ability of the DGS gradient to find a better design. The implementation of the DGS-ES algorithm for TO problem can be summarized as the following steps: 
\vspace{-0.1cm}
\begin{itemize}[leftmargin=13pt]
    \item Make an initial design using Gaussian random noise. 
    
    \item For the given distribution of density, compute the displacement using FEM. 
    
    \item Compute the objective, typically the compliance of this design, and the associated gradient with respect to design changes using DGS-ES. If the change is smaller than the specific threshold, stop the iteration, otherwise continue.
    
    \item Compute the update of the density variable, by solving the MMA approximation subproblem using a dual or primal-dual method. 
    
    \item Repeat the iteration loop. 
\end{itemize}

Our implementation is built on the python implementation published by \cite{andreassen2011efficient} at \url{http://www.topopt.mek.dtu.dk/Apps-and-software/Topology-optimization-codes-written-in-Python} and we choose density-based filtering with filter radius $f_r$=1.5 to avoid the numerical instability like checkerboard problem. For MMA optimizer, we use the python implementation from \url{https://github.com/arjendeetman/TopOpt-MMA-Python} with the default hyper-parameters defined in MMA. To reduce the computational cost, we use mpi4py to parallelly run the TO example in \S 4.2 on multiple cores workstation. 

\subsubsection{Hyper-parameters of each method in TO}
The hyper-parameters of DGS-ES in TO design problem are $M=5$, $\alpha=0.1$, $r=0.25$, $\beta=0.2$ and $\gamma=0.01$. Note that there is no learning rate $\lambda$ in this case because the update step of design variable is achieved by MMA optimizer. The hyper-parameters for other compared methods include: (1) ES-Bpop: $\sigma = 0.25$ and $\lambda_t = 0.99\lambda_{t-1}$ with $\lambda_0 = 0.1$; (2) ASEBO: the population size is $14 \approx 4+3\log(2400)$, $\sigma = 0.1$ and $\lambda_t = 0.99\lambda_{t-1}$ with $\lambda_0 = 0.1$; (3) IPop-CMA: \texttt{restarts}=9, \texttt{restart\_from\_best=False}, \texttt{incpopsize$~ =2$}, \texttt{$\sigma_0=0.25$}; (4) Nesterov: $\lambda_t = 0.99\lambda_{t-1}$ with $\lambda_0 = 0.01$; (5) FD: $\lambda_t = 0.99\lambda_{t-1}$ with $\lambda_0 = 0.01$.

\subsection{Additional discussion}
The results in \S 4.2 example demonstrate that the IPop-CMA underperforms all other methods since it fails to solve such constrained optimization well. The main challenge for IPop-CMA is its limitation in efficiently handling constraints because of its sample-based update mechanism, and in effectively incorporating with other optimizer, such as MMA.  To satisfy the volume constraint in TO problem, IPop-CMA uses Lagrangian penalty as regularization term to enforce the constraints but loses its own capability to seek optimization. The regularization coefficient is highly sensitive to balance the penalty term and loss function in pycma (v3.0.3 available at \url{https://github.com/CMA-ES/pycma}) with constraints. In other words, it is difficult for IPop-CMA to address constrained optimization by using a simple penalty approach. This has also been the challenges of CMA-based approaches in complex real-world engineering applications.

\section{Additional information on the hydrology example in \S 4.3}
We provide detailed information on the implementation of the hydrology example. We consider a 2D square aquifer domain, denoted by $\mathcal{D} = [0,1] \times [0,2]$. The domain the aquifer is discritized into a $100 \times 200$ mesh. The partial differential equation (PDE) governing the groundwater flow is 
\begin{eqnarray*}
-\nabla \cdot (a(\bm x) \nabla u(\bm x)) &  = 0 & \text{ for } \bm x \in (0,1)\times (0,2) \notag \\
\nabla u(\bm x) & = 0 & \text{ for } x_2=0 \text{ and } x_2=2 \\
u(\bm x) & = 10 & \text{ for } x_1=0 \notag \\
u(\bm x) & = 100 & \text{ for } x_1=1, \notag
\end{eqnarray*}
where $u(\bm x)$ is the hydraulic head field and $a(\bm x)$ is the hydraulic conductivity field. The no-flow boundary condition $\nabla u =0$ is imposed to the top $(x_2=2)$ and bottom $(x_2=0)$ boundaries, and the constant hydraulic head condition is applied to the left ($x_1=0$) and right ($x_1=1$) boundary. The goal is to infer $a(\bm x)$ using the sampled data of $u(\bm x)$.

The ground-truth of the hydraulic conductivity field (shown in Figure 7 (Left)) is generated by a sequential Gaussian sampling. As discussed in \S 4.3, we use a full connected neural network (FNN) to approximate the logarithm of $a(\bm x)$, i.e., $\text{FNN}(\bm x; \bm w) \approx \log(a(\bm x))$. The loss function is defined by
\[
Loss = \frac{1}{S}\sum_{s=1}^S\big[ \text{MODFLOW}(\text{FNN}(\bm x))(\bm x_s) - u(\bm x_s)\big]^2,
\]
where $S = 50$, and $u(\bm x_s)$ for $s=1,\ldots, S$ are the hydraulic head data sampled at 50 random locations. The simulator MODFLOW($\cdot$) maps the predicted hydraulic conductivity field $\text{FNN}(\bm x))$ to the hydraulic head by solving the above PDE. 

The hyper-parameters are set as follows. DGS-ES: $M=5, \alpha = 0.1, r = 0.1, \beta = 0.1, \gamma = 0.001$ and $\lambda_t = 0.99\lambda_{t-1}$ with $\lambda_0 = 0.1$; ES-Bpop: $\sigma = 0.1$ and $\lambda_t = 0.99\lambda_{t-1}$ with $\lambda_0 = 0.1$; ASEBO: $\sigma = 0.1$, $\lambda_t = 0.99\lambda_{t-1}$ with $\lambda_0=0.1$; IPop-CMA: \texttt{restarts=9}, \texttt{restart\_from\_best=False}, \texttt{incpopsize=2}, \texttt{$\sigma_0=0.3$}; Nesterov: $\lambda_t = 0.99\lambda_{t-1}$ with $\lambda_0 = 0.05$; FD: $\lambda_t = 0.99\lambda_{t-1}$ with $\lambda_0 = 0.05$.

\section{Additional discussion on asymptotic consistency}
%
We provide additional results to support the discussion on the asymptotic consistency in \S 3.2. Recall that $F \in C^{1,1}(\mathbb{R}^d)$ if there exists $L > 0 $ such that
 $
 \|\nabla F(\bm x + \bm \xi) - \nabla F(\bm x)\| \le L\|\bm \xi\|, \, 
 \forall \bm x ,\, \bm \xi\in \mathbb{R}^d.
 $
Also, recall the error of one-dimensional GH quadrature
\begin{align}
\label{GH_error2}
\hspace{-0.1cm}\big|(\widetilde{\mathscr{D}}^M- \mathscr{D})[G_\sigma] \big| \le C\frac{M\,!\sqrt{\pi}}{2^M(2M)\,!} \sigma^{2M-1}, 
\end{align}
where $C>0$ is a constant independent of $M$ and $\sigma$. 
Given a unit vector $\bm \xi \in \mathbb{R}^d$, we define $\nabla_{\bm \xi} F(\bm x)$ the partial derivatives of $F$ at $\bm x\in \mathbb{R}^d$ in direction $\bm \xi$.
We say $F$ is a strongly convex function if there exists a positive number $\tau$ such that for any $\bm x,\bm \xi\in \mathbb{R}^d$,  
$    F(\bm x + \bm \xi) \ge F(\bm x) + \langle \nabla F(\bm x),\bm \xi \rangle + \frac{\tau}{2} \|\bm \xi\|^2
$. We call $\tau$ the {convexity parameter} of $F$. We prove below the estimate on the difference between DGS estimator $\widetilde{\nabla}^M_{ \sigma, \bm \Xi}[F]$ and $\nabla F$. 

\begin{proposition}
\label{lemma:grad_est}
Let  $\bm \Xi = \{\bm \xi_1, \ldots, \bm \xi_d\}$ be a set of orthonormal vectors in $\mathbb{R}^d$ and $F$ be a function in $C^{1,1}(\mathbb{R}^d)$. Then 
\begin{align}
& \| \widetilde{\nabla}^M_{ \sigma, \bm \Xi}[F](\bm x) -  {\nabla}F(\bm x) \|^2  
 \le \frac{2 C^2 {\pi} d (M\,!)^2}{4^M((2M)\,!)^2}  \sigma^{4M-2} + 32 d L^2  \sigma^2.  \label{lemma:approx_dev0}
\end{align}
\end{proposition}

\textit{Proof.} First, adapting \cite[Lemma 3]{NesterovSpokoiny15} to $1$-dimensional Gaussian smoothing, for any $\bm \xi$ being a unit vector in $\mathbb{R}^d$, there holds
\begin{align}
\left|{ \mathscr{D}}\left[G_{\sigma}(0 \, |\, \bm x, \bm \xi)\right] - \nabla_{\bm \xi} F(\bm x) \right| \le 4\sigma L. \label{lemma:est0}
 \end{align}
From \eqref{GH_error2} and \eqref{lemma:est0}, we have
\begin{gather}
\label{lemma:approx_dev0b}
\begin{aligned}
 \left|\widetilde{\mathscr{D}}^M[G_{\sigma}(0 \, |\, \bm x, \bm \xi_i)] - \nabla_{\bm \xi_i} F(\bm x) \right|^2 
 & \le  2\left|(\widetilde{\mathscr{D}}^M - \mathscr{D})[G_{\sigma}] \right|^2  + 2 \left|{ \mathscr{D}}\left[G_{{\sigma}}\right] - \nabla_{\bm \xi_i} F(\bm x) \right|^2
\\
& \le  \frac{2C^2 (M\,!)^2{\pi}}{4^M((2M)\,!)^2}\, \sigma^{4M-2} + 32\sigma^2 L^2, \ \ \forall i\in \{1,\ldots,d\}.
\end{aligned}
\end{gather}
Summing \eqref{lemma:approx_dev0b} from $i=1$ to $d$ gives \eqref{lemma:approx_dev0}. 
$\square$

Let $N$ be the total number of function evaluations. In our DGS-ES algorithm, $N=Md$. 
An immediate consequence of Proposition \ref{lemma:grad_est} is that given a positive $\varepsilon$, $\sigma \le {\varepsilon}/(4L\sqrt{d})$
and $N  \ge d\log({2d}/{\varepsilon^2})$ are sufficient to obtain $\| \widetilde{\nabla}^M_{ \sigma, \bm \Xi}[F] -  {\nabla}F \| \le \varepsilon$. Now, we compare these with the condition on $N$ such that $\|g(\bm x) - \nabla F(\bm x)\| \le \varepsilon$, where $g(\bm x)$ is an MC-based gradient estimator for $\nabla F_{\sigma}$. For simplicity, we focus on the condition for $\|g(\bm x) - \nabla F_{\sigma}(\bm x)\| <\varepsilon$, which is actually weaker because it does not count for the discrepancy between $\nabla F_{\sigma}(\bm x)$ and $\nabla F(\bm x)$. Let $\text{Var}[g(\bm x)]$ be the variance of $g$, applying Chebyshev inequality, one has
$$
\mathbb{P}\left(\|g(\bm x) - \nabla F_{\sigma}(\bm x)\| > \varepsilon \right) \le \frac{d ~ \text{Var}[g(\bm x)]}{\varepsilon^2}, 
$$
therefore,  
$
 \|g(\bm x) - \nabla F_{\sigma}(\bm x)\| \le \varepsilon
$
with probability exceeding $1-\eta$ given that $\text{Var}[g(\bm x)] < {\eta\varepsilon^2}/{d}$. It can be shown for forward finite-difference MC gradient estimator that 
$
    \text{Var}[g(\bm x)] \simeq {\|\nabla F(\bm x)\|^2}/{N},
$
see \cite{BCCS19}, thus, this basic estimator requires $N = O( \frac{d\|\nabla F(\bm x)\|^2}{\eta \varepsilon^2 })$. It is possible to reduce the number of function evaluations with various variance reduction techniques, such as antithetic sampling, orthogonalization, control variates. In particular, if the variance of MC estimator is $\kappa \text{Var}[g(\bm x)] $ instead of $\text{Var}[g(\bm x)]$ for some $\kappa < 1$, then $N = O( \frac{\kappa d\|\nabla F(\bm x)\|^2}{\eta \varepsilon^2 })$ is sufficient. However, the proven rate of reduction $\kappa$ is either independent or $O(1)$ on $\delta$ and $\varepsilon$ for most variance reduction techniques, see, e.g., \cite{CRSTW18,Tang2019VarianceRF}, in which cases, the theoretical dependence of $N$ on $d$ and $\varepsilon$ cannot be relaxed.    


With the gradient estimate in Proposition \ref{lemma:grad_est}, we proceed to establish an error analysis for DGS-ES in local regime. We show here a convergence rate of our method in optimizing strongly convex functions. 
\begin{proposition}
\label{thm:main}
Let $F$ be a strongly convex function in $C^{1,1}(\mathbb{R}^d)$, $\bm x^*$ be the global minimizer of $F$ and the sequence $\{\bm x_t\}_{t\ge 0}$ be generated by Algorithm 1 with $\lambda = {1}/({8L})$. Then, for any $t\ge 0$,
\begin{align}
        &  F(\bm x_t) - F(\bm x^*) \le\,   
 \frac{1}{2}L \left[\delta_{\sigma} + \left(1-  \frac{\tau}{16 L}\right)^t  ( \|\bm x_0 - \bm x^*\|^2 - \delta_{\sigma}) \right]. \label{thm:strong_convex}
\end{align}
\begin{align}
 \text{Here,}\qquad \   & \delta_{\sigma} =  
   \left(\frac{128}{\tau^2} + \frac{16}{\tau L}  \right)  L^2 d  \sigma^{2} 
  + \, \left(\frac{8}{\tau^2} + \frac{1}{2\tau L}  \right)\frac{ C^2 (M\,!)^2{\pi}d}{4^M((2M)\,!)^2}   \sigma^{2}.  \label{thm:delta}
\end{align} 
\end{proposition}

\textit{Proof.}
First, we derive an upper bound for $\|\widetilde{\nabla}^M_{ \sigma, \bm \Xi}[F](\bm x)\|$. Recall 
 \begin{align}
     & \|\widetilde{\nabla}^M_{ \sigma, \bm \Xi}[F](\bm x)\|^2  =  \sum_{i=1}^d \left|\widetilde{ \mathscr{D}}^M\left[G_{\sigma}(0 \, |\, \bm x, \bm \xi_i)\right] \right|^2 .  \label{lemma:approx_dev:est1}
\end{align}
Each term inside this sum can be bounded as    
\begin{align*}
     \left|\widetilde{ \mathscr{D}}^M\left[G_{\sigma}(0 \, |\, \bm x, \bm \xi_i)\right] \right|^2 & \le 2  \left|{ \mathscr{D}}\left[G_{\sigma}(0 \, |\, \bm x, \bm \xi_i)\right] \right|^2
    +\, 2  \left|\widetilde{ \mathscr{D}}^M\left[G_{\sigma}(0 \, |\, \bm x, \bm \xi_i)\right] - { \mathscr{D}}\left[G_{\sigma}(0 \, |\, \bm x, \bm \xi_i)\right] \right|^2
    \\
    & \overset{\eqref{GH_error2},\eqref{lemma:est0}}{\le} \,  64\sigma^2 L^2 + 4|\nabla_{\bm \xi_i} F(\bm x)|^2  + \frac{2C^2(M\,!)^2{\pi}}{4^M((2M)\,!)^2} \sigma^{4M-2}.
\end{align*}
Plugging this into \eqref{lemma:approx_dev:est1} gives
\begin{align}
\| \widetilde{\nabla}^M_{ \sigma, \bm \Xi}[F](\bm x)\|^2 \le \,   64 d L^2  \sigma^2   
 +\, 4\sum_{i=1}^d |\nabla_{\bm \xi_i} F(\bm x)|^2 + \, \frac{2C^2{\pi} d (M\,!)^2}{4^M((2M)\,!)^2}  \sigma^{4M-2}.  \label{lemma:approx_dev}
\end{align} 
Denote $r_t = \|\bm x_t - \bm x^*\|$. Then
\begin{align}
 r_{t+1}^2  = & \|\bm x_t - \lambda \widetilde{\nabla}^M_{ \sigma, \bm \Xi}[F](\bm x_t) - \bm x^*\|^2 \label{thm:est1}
\\
 = &\, r_{t}^2 - 2\lambda \langle \widetilde{\nabla}^M_{ \sigma, \bm \Xi}[F](\bm x_t), \bm x_t - \bm x^* \rangle + \lambda^2 \| \widetilde{\nabla}^M_{ \sigma, \bm \Xi}[F](\bm x_t)\|^2\notag
\\
 \overset{\eqref{lemma:approx_dev}}{=} &\, r_{t}^2 - 2\lambda \langle \widetilde{\nabla}^M_{ \sigma, \bm \Xi}[F](\bm x_t)  -  {\nabla}F(\bm x_t), \bm x_t - \bm x^* \rangle - 2\lambda \langle {\nabla}F(\bm x_t), \bm x_t - \bm x^* \rangle  \notag
 \\
  \qquad & + \, 4\lambda^2\sum_{i=1}^d |\nabla_{\bm \xi_i} F(\bm x_t)|^2 + \, 64 \lambda^2 L^2 d  \sigma^2    + \, \frac{2C^2 \lambda^2 (M\,!)^2{\pi}d}{4^M((2M)\,!)^2}   \sigma^{4M-2}. \notag
\end{align}
We proceed to bound the right hand side of \eqref{thm:est1}. First, since $F$ is strongly convex, 
\begin{gather}
\label{thm:strongly_convex}
\begin{aligned}
      & - 2\lambda\langle\nabla F(\bm x_t) , \bm x_t - \bm  x^* \rangle \le 2\lambda F(\bm x^*) 
  - 2\lambda F(\bm x_t)  - {\lambda \tau}\|\bm x^* - \bm x_t\|^2. 
\end{aligned}
\end{gather}
On the other hand, 
\begin{gather}
\label{thm:est4}
\begin{aligned}
 & - 2\lambda \langle \widetilde{\nabla}^M_{ \sigma, \bm \Xi}[F](\bm x_t) -  {\nabla}F(\bm x_t), \bm x_t - \bm x^* \rangle 
 \\
 \le\, & \frac{2\lambda}{\tau} \| \widetilde{\nabla}^M_{ \sigma, \bm \Xi}[F](\bm x_t) -  {\nabla}F(\bm x_t)\|^2 + \frac{\lambda\tau}{2} \| \bm x_t - \bm x^* \|^2 
 \\
 \overset{\eqref{lemma:approx_dev0}}{\le}\, &  \frac{4\lambda}{\tau}\cdot \frac{ C^2 {\pi} d (M\,!)^2}{4^M((2M)\,!)^2}   \sigma^{4M-2}  +  \frac{64\lambda}{\tau} L^2 d \sigma^2  + \frac{\lambda\tau}{2} \| \bm x_t - \bm x^* \|^2. 
\end{aligned}
\end{gather}


Applying an estimate for convex, ${C}^{1,1}$-functions, see, e.g., \cite[Theorem 2.1.5]{Nesterov_book_2004}, gives
\begin{align}
  &  4\lambda^2\sum_{i=1}^d |\nabla_{\bm \xi_i} F(\bm x_t)|^2 =   4\lambda^2\|\nabla F(\bm x_t)\|^2
 \le 8\lambda^2 L(F(\bm x_t) - F(\bm x^*)).  \label{thm:est4b}
\end{align}
Combining \eqref{thm:est1}--\eqref{thm:est4b}, there holds 
\begin{gather}
\label{thm:est5}
\begin{aligned}
  r_{t+1}^2 &  \, {\le} \, r_{t}^2     
- (2\lambda -8\lambda^2 L)( F(\bm x_t) -   F(\bm x^*) )  
 \\
 &  + \left(\frac{64\lambda}{\tau} +  64 \lambda^2 \right) L^2 d  \sigma^2  
  + \, \left(\frac{4\lambda}{\tau} + 2\lambda^2\right)\frac{ C^2 (M\,!)^2{\pi}d}{4^M((2M)\,!)^2}   \sigma^{4M-2}. 
\end{aligned}
\end{gather}
Since $F$ is a strongly convex function, for $\lambda = 1/(8L)$ we have that 
\begin{align}
    & - (2\lambda -  8\lambda^2 L) (F(\bm x_t) -   F(\bm x^*)) 
       \le \,   \frac{\tau}{16 L}\|\bm x_t - \bm x^*\|^2. \label{thm:est5b}
\end{align}
Assuming $ \sigma<1$. We derive from \eqref{thm:est5}, \eqref{thm:est5b} and \eqref{thm:delta} that
$
      r_{t+1}^2 - \delta_{\sigma} \le   \left(1-  \frac{\tau}{16 L}\right) (r_{t}^2 - \delta_{\sigma}),
$
which yields 
$$
r_{t}^2 - \delta_{\sigma} \le \left(1-  \frac{\tau}{16 L}\right)^t  (r_{0}^2 - \delta_{\sigma}). 
$$
Note that $F(\bm x_t) - F(\bm x^*) \le \frac{1}{2}L\|\bm x_t - \bm x^*\|^2$, since $f\in C^{1,1}(\mathbb{R}^d)$, we arrive at the conclusion. 
$\square$

It is worth remarking a few things on the above proposition. First, we obtain the global linear rate of convergence with DGS-ES, which is expected for strongly convex functions.
Second, the result allows random perturbation of $\bm \Xi$ as long as $\bm \Xi$ remains orthonormal. Finally, this proposition proves the scalability of our algorithm in the strongly convex setting. In particular, to guarantee $F(\bm x_t) - F(\bm x^*) \le \varepsilon$, from error estimate \eqref{thm:strong_convex}, we need to choose 
\begin{align*}
\sigma \le O\left(\sqrt{\frac{\varepsilon}{d}}\right),\ &\#\text{ function evaluations}=dM\ge O\left(d\log\left(\frac{d}{\varepsilon}\right)\right),
\\
&\#\text{ iterations }= O\left(\log\frac{1}{\varepsilon}\right).
\end{align*}
This indicates that the number of iterations required by our approach is completely independent of the dimension, while the total number of function evaluations is only slightly higher than nonparallelizable random search approach, e.g., \cite{NesterovSpokoiny15,Bergou2019}. 

The above discussion shows that the performance of the DGS-ES method is consistent with the local gradient estimation methods when $\sigma$ is small, which paves the way for further analysis in the nonlocal regime for which the DGS gradient is designed.


\end{document}